\documentclass[12pt,reqno]{amsart}
\usepackage{amssymb,latexsym,amsmath}
\usepackage{amsthm, color}
\usepackage[all]{xy}

\usepackage[pagebackref,colorlinks,
linkcolor=black,
citecolor=blue,urlcolor=blue,
hypertexnames=false]{hyperref}
\newcommand{\Cf}{\mathrm{C}}
\newcommand{\Co}{\mathcal{C}}

\newcommand{\cF}{\mathcal{F}}
\newcommand{\Hi}{\mathcal{H}}

\newcommand{\cK}{\mathcal{K}}
\newcommand{\cL}{\mathcal{L}}
\newcommand{\cM}{{\mathcal{M}}}

\newcommand{\cS}{\mathcal{S}}

\newcommand{\C}{\mathbb{C}}
\newcommand{\R}{\mathbb{R}}
\newcommand{\Z}{\mathbb{Z}}
\newcommand{\N}{\mathbb{N}}
\newcommand{\pf}
{\noindent{\mbox{\textbf{Proof}.\,}} }
\newcommand{\cf}{\textrm{cf.~}}
\newcommand{\eg}{\textrm{e.g.~}}
\newcommand{\ie}{\textrm{i.e., }}
\newcommand{\cst}{\textit{C}*}
\newcommand{\id}{\mathrm{id}}

\theoremstyle{plain}
\newtheorem{thm}{Theorem}[section]
\newtheorem{cor}[thm]{Corollary}
\newtheorem{lem}[thm]{Lemma}
\newtheorem{prop}[thm]{Proposition}
\theoremstyle{definition}
\newtheorem{definition}[thm]{Definition}

\newtheorem{rem}[thm]{Remark}
\newtheorem{rems}[thm]{Remarks}

\begin{document} 

\title[Strong pure infiniteness of crossed products]
{Strong pure infiniteness of crossed products}

\author[Kirchberg and Sierakowski]{E.~Kirchberg 
and A.~Sierakowski}

\date{\today} 


\begin{abstract} Consider an exact action 
of discrete group $G$ 
on a separable \cst-algebra $A$. 
It is shown that the reduced crossed 
product $A\rtimes_{\sigma, \lambda}  G$ 
is strongly purely 
infinite -- provided that the action of $G$ on any 
quotient $A/I$ by a $G$-invariant closed ideal 
$I\neq A$ 
is element-wise properly outer 
and that the action of $G$
on $A$ is \emph{$G$-separating} 
(\cf Definition\ \ref{def:G.separating}).
This is the first non-trivial sufficient criterion for 
strong pure infiniteness of reduced crossed
products
of 
\cst-algebras $A$ that are not $G$-simple.
In the case 
$A=\mathrm{C}_0(X)$ the notion of a 
$G$-separating 
action corresponds to the 
property that two compact 
sets $C_1$ and $C_2$, 
that  are contained in open subsets
$C_j\subseteq  U_j \subseteq  X$, 
can be mapped by elements 
of $g_j\in G$ 
onto 
disjoint sets  
$\sigma_{g_j}(C_j)\subseteq  U_j$, but we
do not require that 
$\sigma_{g_j}(U_j)\subseteq U_j$.
A generalization
of strong boundary actions
\cite{LacaSpi:purelyinf}  on compact spaces
to  non-unital and non-commutative 
\cst-algebras  $A$
(\cf Definition \ref{def:n.majorizing.action})
is also introduced. It
is
stronger than the notion of $G$-separating 
actions by Proposition 
\ref{prop:1.majorizing.implies.G.separating},
because $G$-separation does not imply 
$G$-simplicity
and there are examples of 
$G$-separating actions
with reduced crossed products that are
stably projection-less and non-simple.
\end{abstract}

\maketitle

\tableofcontents

\section{Introduction}\label{sec:intro}
In this paper we pursue the study of \cst-dynamical 
systems with applications in classification via
equivariant $KK$-theory.
It was shown by the first named author 
that for any two separable nuclear strongly 
purely infinite \cst-algebras, both with primitive 
ideal space isomorphic to a $T_0$-space $X$, the 
algebras are isomorphic if and only if they are 
$KK_X$-equivalent. It is however far from 
understood when \cst-algebra crossed products
$A\rtimes_{\sigma, \lambda}  G$
associated to
\cst-dynamical systems are strongly 
purely infinite in terms of properties 
of the action $\sigma$, 
in particular in the non-simple
case. Our main focus of this work is such
characterisation for crossed products 
that are either 
simple or more generally contain ideals 
coming from arbitrary $G$-invariant ideals 
of the algebra $A$
on which the group $G$ acts.

We begin (in Section \ref{sec.2}) by introducing 
crossed products and by proving the
notation used throughout the paper.

In Section \ref{sec.3} we look at 
results related to the 
ideal structure of crossed products. 
It was shown in 
\cite{Siera2010} that 
{residually} properly
outer 
(Definition 
\ref{def:G.separating})
and exact 
(Definition 
\ref{def:exact})
actions 
$\sigma\colon G\to \mathrm{Aut}(A)$ 
on a separable 
\cst-algebra $A$  have the property
that the lattice of (closed) ideals of 
the reduced crossed product  
$A\rtimes_{\sigma, \lambda}  G$
is naturally isomorphic to the lattice of 
$G$-invariant
ideals of $A$ (by the map $I\mapsto A\cap I$).
We refine this result by showing
that for any exact and residually 
properly outer action $\sigma$ of a discrete group 
$G$ on a separable or commutative \cst-algebra $A$ 
the set $A_+$ is a filling family 
(Definition 
\ref{def:filling.mdiag.family})
for 
$A \rtimes_{\sigma, \lambda}   G$
(which implies that $I\mapsto A\cap I$ is injective, 
see Remark \ref{rems:filing.familiy} for details).

In Section \ref{sec:spi.crossed} we introduce 
the notion of $G$-separating actions
(Definition 
\ref{def:G.separating}). We
show in Theorem \ref{cor:A-filling} that for 
any exact and residually 
properly outer action $\sigma$ of a discrete group 
$G$ on a separable or commutative \cst-algebra $A$
and for any filling family $\cF\subseteq A_+$, the 
crossed product $A \rtimes_{\sigma, \lambda}   G$
is strongly purely infinite if and 
only if $\cF$ has the diagonalization property 
(Definition 
\ref{def:mdiag.family})
in $A \rtimes_{\sigma, \lambda}   G$. 
Applying the work \cite{KirSie2} we obtain
(in Proposition \ref{prop:infinite}) an equivalent 
characterisation of $G$-separating 
actions, from which we can deduce that
$A_+$ has the diagonalization property whenever 
the action on $A$ is $G$-separating. 
By evoking \cite{KirSie2} once again we prove 
our main result:
\begin{thm}\label{thm:main}
Suppose that 
$(A,G,\sigma)$ is a \cst-dynamical 
system, where $G$ is discrete and $A$ is
separable or commutative.

If the action $\sigma$ of $G$ on $A$ 
is exact (Def.~\ref{def:exact}), 
residually properly outer
(Def.~\ref{def:properly.outer})
and  $G$-separating 
(Def.~\ref{def:G.separating}),
then $A\rtimes_{\sigma, \lambda}  G$ 
is strongly purely infinite.
\end{thm}

In Section \ref{sec:abelian.case} we look at 
actions on commutative \cst-algebras. Here we 
characterise the notion of $G$-separating action 
purely in terms of the underlying geometry. More 
specifically we consider actions $\alpha$ of a
discrete group $G$ on a locally compact Hausdorff
space $X$, and denote by $\sigma$
the induced action on $A:=\Cf_0(X)$.
We show (in Lemma \ref{lem:G-separating}) that 
the action of $G$ on $A$ is $G$-separating if 
and only if the following holds:
For every open 
$U_1,U_2\subseteq X$  and 
compact $K_1,K_2\subseteq X$ with 
$K_1\subseteq U_1$, $K_2\subseteq U_2$,
there exist $g_1,g_2\in G$ such that 
$\alpha_{g_1}(K_1)\subseteq U_1$, 
$\alpha_{g_2}(K_2)\subseteq U_2$,
$\alpha_{g_1}(K_1)\cap \alpha_{g_2}(K_2)
 =\emptyset$. This result is what motives
the choice of our terminology 
``$G$-separating''. As a consequence
of this characterisation we obtain in 
Corollary \ref{cor:Abelian}
a characterisation of when a crossed product
$\Cf _0(X)\rtimes_{\sigma, \lambda}  G$
is a strongly purely infinite \cst-algebra
in terms of condition on $\alpha$.

In the final Section 
\ref{sec:general.strong.boundary.actions} 
we consider actions that produce
simple and strongly purely infinite crossed
products. We introduce 
(Definition 
\ref{def:n.covering.action} and 
\ref{def:n.majorizing.action})
the notion
of $n$-majorizing ($n\geq 1$) and $n$-covering 
actions ($n\geq 2$), 
the later for actions on unital \cst-algebras. 
These two notions aim to refine 
results on simple purely infinite crossed 
products in 
\cite{LacaSpi:purelyinf, 
JolissaintRobertson}
where the notion of strong boundary actions 
(Definition 
\ref{def:strong.boundary.action})
and $n$-fillings
actions 
(Definition 
\ref{def:n.filling.action})
was introduced.
We prove in Remark 
\ref{rem:unital.abel.2.cover=1.major=sba} 
our notions are weaker: 
Any $n$-filling action on a unital
\cst-algebra $A$ is $n$-covering, 
and for any action $\alpha$
on a compact spaces $X$ with
more than two points (on which 
strong boundary actions are defined)
the action $\alpha$ is
a strong boundary action if and only
if its adjoint action $\sigma$ on 
$\Cf(X)$ is $1$-majorizing.
Both our notions
are $G$-simple.
\emph{Therefore we call the 
$1$-majorizing actions on not-necessarily
unital or commutative \cst-algebras also 
strong boundary actions}. Despite our 
weaker assumptions we are able to 
prove:
\begin{thm}
\label{thm:n.majorizing.n.covering}
Suppose that the   
\cst-dynamical system  $(A,G,\sigma)$ with
discrete $G$ is  
$n$-majorizing 
(Def.~\ref{def:n.majorizing.action})
for some $n\ge 1$
or $n$-covering
(Def.~\ref{def:n.covering.action})
for some $n\ge 2$, the latter if $A$ 
is unital.
If the action $\sigma$ is 
element-wise  properly outer
(Def.~\ref{def:properly.outer}), 
and $A$ is separable or commutative,
then $A\rtimes_{\sigma,\lambda} G$ 
is simple strongly purely infinite.
\end{thm}

In Section 
\ref{sec:general.strong.boundary.actions} 
we also look at how 
the different
properties relate
to each other. 
In Lemma 
\ref{lem:deny.socle.for.major.cover}
we show that
each $n$-covering action
(for $n\geq 2$)
on a unital \cst-algebra $A$
action is $n$-majorizings,
and the latter properly
(for $n\geq 1$)
implies that the action is 
$(n + 1)$-covering.
In 
Proposition 
\ref{prop:1.majorizing.implies.G.separating} 
we prove that
that any $1$-majorizing action
on a non-unital \cst-algebra $A$
is automatically $G$-separating. 
In Remark 
\ref{rem:unital.abel.2.cover=1.major=sba} 
we prove that
any action on a unital
commutative \cst-algebra $A$ is 
$n$-filling if and only if it is
$n$-covering, and for $n=2$ 
this is again equivalent
to a strong boundary
(\ie$1$-majorizing) action.

We end with a number of remarks, including
a proof of the fact that our notions 
of $G$-separating, $n$-majorizing and
$n$-covering actions can be expressed in terms
of projections when $A$ has real rank zero
(see Remark \ref{rem:case.RR=0}).

We hope that the study 
of crossed products --
even those for actions 
of amenable discrete groups on locally 
compact Polish spaces --
can help to detect possible differences 
between strong and weak pure infiniteness. 
This paper is a very first step in this 
direction, and gives a sufficient 
criterium by conditions
on the action that implies
strong pure infiniteness of 
reduced crossed products.

\section{Preliminaries}\label{sec.2}
We let $A_+$ denote the set of positive 
elements in a \cst-algebra $A$.
We denote 
the positive and the negative part 
of a selfadjoint element $a\in A\,$ by
$a_+:=(|a|+a)/2\in A_+$ and
$a_-:=(|a|-a)/2)\in A_+\,$, where $|a|:=(a^*a)^{1/2}$. 
If $a\in A_+$, 
then $(a-\varepsilon)_+$, the positive part of
$a-\varepsilon 1$ in $\cM (A)$, 
is again in $A_+$. Here  $\cM (A)$ is the 
\emph{multiplier algebra} of $A$. 
This notation
will be used also for functions 
$f\colon \R\to \R$, then \eg 
$(f-\varepsilon)_+(\xi)= \max (f(\xi)\, 
-\varepsilon\, , 0)\,$.
A subset $\cF\subseteq  A_+$ is 
invariant under \emph{$\varepsilon$-cut-downs} 
if for each $a\in \cF$ and 
$\varepsilon \in (0,\| a\|)$
we have
$(a-\varepsilon)_+\in \cF$.
The minimal unitalisation of $A$
is denoted $\tilde{A}$.
Restriction of a map $f$ to $X$
is denoted $f|X$.
We
let $\Cf_c (0,\infty]_+$ 
denote the set of all non-negative continuous
functions $\varphi$ on $[0,\infty)$ with 
$\varphi | [0,\eta] = 0$ 
for some $\eta \in (0,\infty)$, 
such that $\lim_{t\to \infty} \varphi(t)$ exists.

\begin{rems}\label{rems:b-eps=d*ad} 
\noindent 
(i)\,
Suppose that $\,a,b\in A_+\,$ and 
$\,\varepsilon > 0\,$ satisfy  
$\| \,  a-b\, \|  \, <\varepsilon\,$. 
Then the positive part 
$\,(b-\varepsilon)_+\in A\,$ of 
$\,(b-\varepsilon \cdot 1)\in\cM (A)\,$ 
can be decomposed into 
$\,d^*ad=(b-\varepsilon)_+\,$ with 
some contraction $d\in A\,$
(\cite[lem.~2.2]{KirRorOinf}).
\smallskip

\noindent 
(ii)\,
Let $\tau \in [0,\infty)$ and 
$0\leq b\leq a+\tau \cdot 1$ 
(in $\cM  (A)$), 
then for every 
$\varepsilon >\tau$ there is a 
contraction
$f\in A$ such that $(b-\varepsilon)_+=f^*a_+f$.
(See
\cite[lem.~2.2]{KirRorOinf} and
\cite[sec.~2.7]{BlanKir2}.)
\end{rems}

We abbreviate
\cst-dynamical systems by $(A,G,\sigma)$
with discrete groups $G$.
We denote by $e$ the unit of $G$,
and consider only 
\emph{closed and two-sided} ideals of $A$.
The reduced (resp.\ the full)
crossed product associated to 
$(A,G,\sigma)$ is denoted by 
$A \rtimes_{\sigma,\lambda}  G$ 
(resp. $A \rtimes_\sigma G$).
The norm 
on $A \rtimes_{\sigma,\lambda}  G$
will be sometimes written as 
$\|\, \cdot \, \|_\lambda\,$
if it is necessary to distinguish it from other 
norms. Let $\mathcal{I}(A)$ 
denote the lattice of  ideals in a 
\cst-algebra $A$.
The map
$ \eta\colon A\to A\rtimes_{\sigma} G$
means the natural embedding into 
the full crossed product.
Let 
$\pi_\lambda\colon A\rtimes_{\sigma} G
\to  A \rtimes _{\sigma, \lambda}  G$
be the natural epimorphism.
We will sometimes 
suppress the canonical inclusion maps 
$ \eta\colon A\to A\rtimes_{\sigma} G$  
and
$\pi_\lambda\circ \eta\colon A \subseteq 
A \rtimes _{\sigma, \lambda}  G$.
Let $U$ denote the 
canonical  unitary representation
$
U\colon G\to \mathcal{M}(A\rtimes_{\sigma} G)
$.
Notice here that the linear span of 
$\eta(A)U(G)$ is
is a dense *-subalgebra of  
$A\rtimes_{\sigma} G$.
We denote by 
$U_{\lambda}\colon G \to 
\mathcal{M}(A\rtimes_{\sigma, \lambda} G)$
the regular representation 
for some more precise explanations.
The same happens with 
$\eta_\lambda := \pi_\lambda \circ \eta$.

The set  $\Cf_c(G,A)$ consists of 
the maps  $f\colon G\to A$
with finite support $F:=G\backslash  f^{-1}(0)$.
There is a natural linear
embedding of $\Cf _c(G,A)$ into 
$A\rtimes_{\sigma} G\,$
by  canonical identification of $f\colon G\to A$ 
(of finite support)
with an element of $A\rtimes_{\sigma} G\,$: 
Let
$F\subseteq G$ be a finite subset, 
with $f(g)=0$ for $g\not\in F$.
Then $f$ will be identified with the element 
$\sum_{g\in F} \eta(a_g) U(g)$
of $A\rtimes_{\sigma} G\,$, where $a_g:=f(g)$.
In this way $\Cf _c(G,A)$ 
becomes a *-subalgebra of 
$A\rtimes_{\sigma} G\,$ that contains $A$.
The natural \cst-morphism 
$\pi_\lambda \colon A\rtimes_{\sigma} G\to 
A \rtimes _{\sigma, \lambda}  G\,$
is faithful on $\Cf _c(G,A)$,  
and we do not distinguish 
between 
$$
\pi_\lambda( \sum_{g\in F} \eta(a_g) U(g))=
\sum_{g\in F} \eta_\lambda (a_g) U_\lambda(g)
\,$$
and $\sum_{g\in F} \eta(a_g) U(g)$. 
In particular, 
$\eta(a)\in A\rtimes_{\sigma} G$ and 
$
\eta_\lambda(a)\in 
A \rtimes _{\sigma, \lambda} G
$
will be denoted simply by $a\in A$,
and $U_\lambda(g)$ might be 
denoted $U(g)$.

We now recall the  
conditional expectation
$E \colon A\rtimes_{\sigma} G 
\to  \eta(A) \cong A\,$:
The map 
$E_{\mathrm{alg}}\colon \, \Cf _c(G,A) \to  A\,$,\, 
$\,\sum_{g\in F} a_g U(g) \mapsto a_e \,$,
extends by continuity to a 
\emph{faithful conditional expectation} 
$E_\lambda  \colon 
A \rtimes _{\sigma, \lambda}  G \to A$. 
In particular $E_\lambda $ 
is a completely positive contraction, 
$
E_\lambda (A \rtimes _{\sigma, \lambda}  G)=
A
$,
and $E_\lambda (b)=0$ imply $b=0$ for 
$b\in  (A \rtimes _{\sigma, \lambda}  G)_+$.
Since $A$ is also contained 
in its full crossed product
$A\rtimes _{\sigma} G$, 
we can use the natural
epimorphism  
$A\rtimes _{\sigma} G\to 
A\rtimes _{\sigma, \lambda}  G$
to define $E$ by 
$E:= E_\lambda \circ \pi_\lambda$ as a (not 
necessarily faithful) conditional expectation 
from $A\rtimes _{\sigma} G$ onto its 
\cst-subalgebra $A$.

\section{Proper outerness and ideal structure}
\label{sec.3}
In this section we look at conditions on a 
\cst-dynamical system $(A,G,\sigma)$ ensuring that 
the set $A_+$ is a filling family for 
$A \rtimes_{\sigma, \lambda} G$ in the sense of 
Definition 
\ref{def:filling.mdiag.family}.
This in particular implies that 
the is a on-to-one correspondence between ideals 
of $A\rtimes_{\sigma, \lambda}  G$ and 
$G$-invariant ideals of $A$, but (as wee shall 
see) also applies to the verification of when 
a crossed product is strongly purely infinite. 
Proper outerness of the automorphisms $\sigma_t$ 
of $A$
defining the action $\sigma$ turns out also to 
be a crucial ingredient.  We recall the 
definition below.
\begin{definition}
\label{def:properly.outer} 
Suppose that $(A,G,\sigma)$ is a 
\cst-dynamical system and that $G$ is discrete. 
The action $\sigma$ will be called 
\emph{element-wise properly outer}  if, for each 
$g\in G\backslash \{ e\}$, the automorphism 
$\sigma_g$ of $A$ is properly outer in the sense 
of \cite[def.~2.1]{Elliott.prop.outer}, \ie  
$\|\,\sigma_g | I \,-\,\mathrm{Ad}(U)\,\|\,=2$ 
for any  $\sigma_g$-invariant non-zero ideal 
$I$ of $A$ and any unitary $U$ in the multiplier 
algebra $\cM(I)$ of $I$. See also 
\cite[thm.~6.6(ii)]{OlePed3}.

We call  here an action $\sigma$ 
\emph{residually properly outer} if, for every 
$G$-invariant ideal $J\not=A$ of $A$, the 
induced action $[\sigma]_J$ of $G$ on $A/J$ is 
\emph{element-wise} properly outer.
\end{definition}

\begin{rems}\label{properly.outer}
(i)\,
Notice that  \emph{element-wise} 
proper outerness passes to subgroups,
i.e., for
each subgroup $H\subseteq  G$
the system $(A,H, \sigma | H)$ is
element-wise properly  outer on $A$ if 
$(A, G, \sigma)$ is element-wise properly outer.
But \emph{residual} proper outerness 
does not necessarily pass to subgroups. 
The system
$(A, H, \sigma | H)$ is not necessarily residually
properly outer, if $(A,G,\sigma)$ is
residual proper outer, because possibly
there could  be
more $H$-invariant
ideals than $G$-invariant ideals of $A$.
\smallskip

\noindent 
(ii)\,
If $A$ is non-commutative,
then topological freeness of
$(A,G,\sigma)$
in sense of  \cite[def.~1]{ArchSpiel} is --
at least formally -- stronger than the
assumption of element-wise
proper outerness of $(A,G,\sigma)$ in
Definition \ref{def:properly.outer}
(\cf  \cite[prop.~1]{ArchSpiel}).
We do not know examples where
they actually differ.
Thus, for non-commutative $A$,
``essential freeness''
of the corresponding action of
$G$ on $\widehat{A}$
in the sense of
\cite[def.~1.17]{Siera2010}
(inspired by
 \cite[def.~4.8]{Ren:fixed})
is  -- formally -- stronger
than our residual proper outerness of
$(A,G,\sigma)$.
\smallskip

\noindent 
(iii)\,
If $G$ is countable and acting on 
$\Cf_0(X)$, 
one can show 
-- using the Baire property
of $X$ --  that  elementwise proper outerness is
 the same as the 
requirement (for the action 
$\alpha$ of $G$ on $X$ with 
$\sigma_g(f):= f\circ \alpha_{g^{-1}}$)
that points with trivial fix-point subgroup
(trivial isotropy) 
are dense in $X$, \ie
Definition 
\cite[def.~1.17]{Siera2010}
holds. 
We can reformulate this as:
stability subgroups of 
non-empty \emph{open} subsets
are trivial.
\end{rems}

\begin{rem}
\label{rem:OlePed.Abelian}
We recall \cite[lem.~7.1.]{OlePed3} 
(\cf also 
\cite[lem.~3.2]{Kishimoto.outer.auto}):\\
{\it 
If $\alpha_1, \alpha_2,..., \alpha_n$
are properly outer automorphisms of a 
separable \cst-algebra 
$A$, there is, for each  
$a_0, a_1, a_2 ,..., a_n \in \widetilde{A}$, 
with $0\not= a_0 \ge 0$, and each 
$\varepsilon > 0$, an element
$x\in A_+$  with $\| x \|= 1$ such that
}
$$
\| xa_0x\| > \| a_0\| -\varepsilon\,, \quad 
\| xa_i \alpha_j(x)\| < \varepsilon\,, \quad
1\leq i,j\leq n\,.
$$

If $A$ is \emph{commutative}, \ie 
$A\cong \mathrm{C}_0(X)$
for $X=\widehat{A}\subseteq A^*$, 
then it is not necessary to suppose 
that $A$ is separable
in the quoted  lemma of 
D.\ Olesen and G.\ Pedersen
(Compare also \cite{Exel}):
An automorphism 
$\sigma \in \mathrm{Aut}(A)$ 
is properly outer, 
if and only if,
for \emph{every} open subset 
$\emptyset\not= U\subseteq  X$ 
there exists 
$y\in U$ with 
$\widehat{\sigma}(y)\not=y$. 
Thus, for every finite set 
$S\subseteq  \mathrm{Aut}(A)$ 
of properly outer automorphisms, 
every non-empty open subset 
$U\subseteq  X$ contains 
a non-empty open subset 
$V\subseteq  U$ with 
$\widehat{\sigma}(V)\cap V=\emptyset$
for all $\sigma \in S$.
If one takes 
$U:=a_0^{-1} (\| a_0\,\| -\varepsilon, \infty)$
and  non-empty 
$V\subseteq  U$ as above, then
each $x\in \Cf_0(X)$ with 
$\| x\|=1$ and 
support in $V$ satisfies
$\| xa_0x\| >  \| a_0\| -\varepsilon$ and  
$x\sigma(x)=0$ for 
$\sigma \in S$.
\end{rem}

\smallskip 

The following 
Lemma \ref{lem:SharpOlePed3}
is a suitable  modification of proofs
of \cite[lem.7.1,thm.7.2]{OlePed3}. 
It has been proved in 
\cite{ArchSpiel} 
under the stronger assumption  
that the action 
$\sigma$ is topologically free, and 
part (iii)
has been shown in 
\cite[thm.~4.1]{KawamuraTomiyama}
even to be equivalent to the 
topological freeness of the action
if $A$ is commutative and unital  
and $G$ is amenable.
Compare also 
Remark \ref{rem:converse.direction} for a
``residual'' version.

\begin{lem}\label{lem:SharpOlePed3}
Suppose that $A$ 
is separable or commutative, 
and that the action 
of $G$ on $A$ 
is element-wise properly outer.
\begin{itemize}
\item[(i)] For every 
$b\in (A\rtimes_{\sigma} G)_+$ 
with $E(b)\not=0$ and $\varepsilon>0$
there exist $x\in A_+$ satisfying that 
$$
\|x\|=1\,, \; \; \; 
\|xbx-xE(b)x\|<\varepsilon \,, \; \; \;  
\|xE(b)x\|>\|E(b)\|-\varepsilon
\,.$$
This holds also for 
$b\in (A \rtimes _{\sigma, \lambda}  G)_+$
and $E_\lambda$ in place of $E$.

\item[(ii)] If 
$h\colon A\rtimes _{\sigma} G\to \cL (\Hi )$
is a *-representation
such that  $h | A$ is faithful,
then $\| h(b) \| \ge \| E(b) \|$ for all 
$b\in (A\rtimes _{\sigma} G)_+$.

In particular,\  the kernel of
$h$ is contained in the kernel $I_\lambda$
of the natural epimorphism 
$\pi_\lambda \colon 
A\rtimes _\sigma  G\to 
A\rtimes _{\sigma, \lambda}  G$.

\item[(iii)]
Every non-zero ideal of 
$A\rtimes _{\sigma, \lambda}  G$,
has non-zero intersection with $A$.
\end{itemize}
\end{lem}
\pf  
(i):
Let 
$b\in (A\rtimes_{\sigma} G)_+$ 
with $E(b)\not=0$,
and  $\varepsilon>0$.  

Let $a_0:=E(b)$.
Since $\Cf _c(G,A)$ is dense in 
$A \rtimes _{\sigma} G$,
there exists 
$a'=c_0+\sum_{j=1}^m a_jU( g_j)
\in \Cf _c(G,A)$  
with $g_i^{-1}g_j\not=e$ and 
$g_j\not= e$ for 
$i \not= j\in\{1,\ldots,m\}$,
and $\| a'-b \|  <\varepsilon/6$.
Since $E$ is a contraction, it follows that 
$\|b-a\|<\varepsilon/3$  and  
$E(a)=a_0=E(b)$
for $g_0:=e$ and
$a:= a_0+ (a'-c_0) =\sum_{j=0}^m a_jU(g_j)$.
By  \cite[lem.~7.1]{OlePed3} and  
Remark \ref{rem:OlePed.Abelian}
 there  exists 
$x\in A_+$ with
$\|x\|=1$,  
$\|xE(a)x\|>\|E(a)\|-\varepsilon/3m$, and 
$\|xa_j\sigma_{g_j}(x)\|<\varepsilon/3m$ 
for $g_j\neq e$, $j=1,\ldots ,m\,$. 
In particular, 
$$
\| xE(b)x \| =\| xE(a)x \| >  
\| E(a)\| -\varepsilon =\| E(b) \| -\varepsilon
\,.$$
Since  
$\|xa_jU(g_j)x \|= \|xa_j\sigma_{g_j}(x)\|\,$ 
we get in $A\rtimes _\sigma G$  that 
$$
\|x(a-E(a))x\|
\leq 
\sum_{g_j\neq e} 
\|xa_j\sigma_{g_j}(x)\|\leq \varepsilon/3
\,.
$$
Thus, in the full crossed product 
$A\rtimes _\sigma G$ we have
$$\|xbx-xE(b)x\| \leq 
\|x(b-a)x\|+\|x(a-E(a))x\|+\|x(E(a)-E(b))x\| 
<\varepsilon\,.
$$ 
The same arguments work 
for $b\in (A \rtimes _{\sigma, \lambda}  G)_+$
and $E_\lambda$ in place of $E$.

(ii): 
The restriction of $h$ to 
$A\subseteq A\rtimes _{\sigma} G$ is faithful,
hence $\| h(a)\| =\| a\|$ for all $a\in A$.
Let $b\in (A\rtimes _{\sigma} G)_+ $ be given. 
If $E(b)=0$ then $\| h(b)\| \ge \| E(b)\|$.
If $ E(b)\not=0$, then select 
$x\in A_+$ as in (i). It follows
that  
$\| h(xE(b)x)\|=\| x E(b)x\| \ge 
\| E(b)\|-\varepsilon\,$.
On the other hand, 
$\| h(b)\| \ge \| h(x)h(b)h(x)\|=\| h(xbx)\|$ and  
$\varepsilon > \| xbx-xE(b)x\| \ge 
\| h(xbx)- h(xE(b)x)\|$.
Thus $\|h(b)\| +\varepsilon \ge \| h(xE(b)x)\|$, 
and 
$\| h(b)\|+ 2\varepsilon \ge \| E(b)\|$ 
for all $\varepsilon >0$.

Since $E= E_\lambda \circ \pi_\lambda$,
we have $b\in (A\rtimes _{\sigma} G)_+$
and $E(b)=0$ implies that $b$ is contained in the 
kernel of 
$\pi_\lambda \colon A\rtimes _{\sigma} G \to 
A\rtimes _{\sigma, \lambda}  G\,$.
In particular, if 
$h\colon A\rtimes _{\sigma} G\to \cL (\Hi )$
is any *-representation with 
$\| h(b) \| \ge \| E(b)\|$
for all $b\in (A\rtimes _{\sigma} G)_+$, 
then the kernel
$h^{-1}(0)$ of $h$  
is contained in the kernel of the
natural epimorphism
$
\pi_\lambda \colon A\rtimes _{\sigma} G \to 
A\rtimes _{\sigma, \lambda}  G
$.

(iii):
Let $I$ a closed ideal of 
$A\rtimes _{\sigma, \lambda}  G$ 
with $I\cap A=\{ 0\}$, consider 
$(A\rtimes _{\sigma, \lambda}  G)/I$
as a \cst-subalgebra of some $ \cL (\Hi )$, 
and let 
$h\colon A\rtimes _{\sigma} G\to \cL (\Hi )$
the corresponding representation with kernel
$h^{-1}(0)= 
J:=\pi_\lambda^{-1}(I)\supseteq\pi_\lambda^{-1}(0)$. 
Then $h$ is faithful
on $A$ and, therefore, satisfies 
$\pi_\lambda^{-1}(0)\supseteq h^{-1}(0)$.
It follows $I=\pi_\lambda(h^{-1}(0))=\{ 0\}$.
\qed
\begin{definition}[{\cite[def.~1.2]{Siera2010}}]
\label{def:exact}
Suppose that $(A,G,\sigma)$ 
is a \cst-dynamical system with locally compact $G$. 
The action $\sigma$  of $G$ on $A$ is \emph{exact},  
if,  for every $G$-invariant ideal $J$ in $A$, 
the sequence 
$0\to\,  J\rtimes_{\sigma | J, \lambda} G \,\to\, 
A\rtimes_{\sigma, \lambda} G\, \to\, 
A/J \rtimes_{[\sigma]_J, \lambda} G\,\to 0$ is 
short-exact.
\end{definition}

\begin{rems}\label{exact}
(i)
Recall that a locally compact group $G$ is 
\emph{exact}, if and only if,
every action 
$\sigma \colon G\to \mathrm{Aut}(A)$
is exact. If $G$ is \emph{discrete},  
then this is equivalent to
the exactness of the \cst-algebra 
$C^*_\lambda (G)$,   
\cf \cite{KirWa.Exact.Groups.Cont.Bund}.
This applies to all amenable groups 
$G$, \eg $G=\Z$.
Under Definition \ref{def:exact}
each minimal (= $G$-simple) action is exact.
In particular, non-exact discrete
groups can have exact (and faithful) actions.
\smallskip

\noindent 
(ii)
Let $F$ denote the (small) 
\emph{Thompson group}
and 
$\rho\colon F\to \mathrm{Homeo}(\R)$
the minimal action of $F$ 
(or only of its commutator subgroup
$F'$) on the real line $\R$ 
as described by
Haagerup and Picioroaga 
in \cite[rem.~2.5.]{HaagerupPicioroaga}.
One can use $\rho$ 
to construct
a $F$-separating, non-minimal 
and exact action $\alpha$
of $F$ (or $F'$) on the 
disjoint union of two 
lines $X:=\R \cup (i+\R)\subseteq \C$ 
if one considers the restriction
$\alpha(g):=\beta(g)|X$
to $X$ of the action 
$g\in F\to \beta(g)$ on
$\C$ given by
$\beta(g)(s+it):= \rho(g)(s)+it$
(\,\footnote{\,
This action is not topologically 
free.}\,).
It is at present  unknown whether  
the Thompson group $F$ 
is exact or not, \cf 
\cite{Arzhantseva.Guba.Sapir,
FarleyDS.Thompson,
Guentner.Kaminker}.
\smallskip

\noindent 
(iii)
It is not known if Gromov's examples 
of non-exact groups can
have non-exact actions on 
\emph{commutative} \cst-algebras.
It is likely that it has to do with still 
missing non-trivial 
\emph{geometric}
conditions for $G$-actions on 
locally compact  spaces $X$
that are equivalent to the 
\emph{exactness} of the adjoint action 
$\sigma\colon G\to \mathrm{Aut}(\Cf_0(X))$
given by $\sigma_g(f):= 
f\circ \alpha_{g^{-1}}$.
Therefore we use the trivial and 
non-geometric definition and define 
$\alpha$ to be exact 
if
its adjoint action $\sigma$ on $\Cf_0(X)$ is exact.
\end{rems}

\begin{rem}\label{rem:from.Sira}
Combination of Lemma 
\ref{lem:SharpOlePed3}(iii)
and of the exactness of an
action 
$\sigma\colon G\to \mathrm{Aut}(A)$ 
on a separable 
or commutative \cst-algebra $A$ shows
that the lattice of (closed) ideals of 
the reduced crossed product  
$A\rtimes_{\sigma, \lambda}  G$
is naturally isomorphic to the lattice of 
$G$-invariant
ideals of $A$ (by the map $J\mapsto A\cap J$),
if $\sigma$ is 
\emph{exact and residual properly outer}.
(See \cite[Remark 2.23]{Siera2010} 
for details.)
\end{rem}

\begin{thm}\label{thm:A-filling}
Let $(A,G,\sigma)$  a \cst-dynamical system,
with discrete  
$G$ and separable or commutative $A$.  
If  the action $\sigma$ of $G$ on $A$
is exact and residually properly outer,
then the elements of $A_+$ build a \emph{filling 
family} for $A \rtimes_{\sigma, \lambda}   G$
in the sense of Definition 
\ref{def:filling.mdiag.family}.
\end{thm}

\pf
We show that for every hereditary 
\cst-subalgebra
$D$ of $A\rtimes_{\sigma, \lambda}  G$ 
and every (closed) ideal
$I $ of $A\rtimes_{\sigma, \lambda}  G$ with 
$D\not\subseteq I $ there
exist $f\in A_+ $ and 
$z\in A\rtimes_{\sigma, \lambda}  G$
such that $z^*z\in D$ and  $zz^*=f \not\in I $.

Suppose that $D$ is  a hereditary 
\cst--subalgebra of 
$A\rtimes_{\sigma, \lambda}  G$ 
and that 
$I $ is an  ideal of 
$A\rtimes_{\sigma, \lambda}  G$ 
with $D \not\subseteq I $.
Let  $J:= I \cap A$, then $J$ is a 
$G$-invariant ideal of $A$ with 
$J\rtimes _{\sigma | J, \lambda} G\subseteq I $  
and 
$g\in G \mapsto [\sigma_g]_J$
is an \emph{element-wise properly outer} 
action on $A/J$
by our assumptions on $\sigma$.
We denote 
this action
by $\alpha$, \ie
$\alpha_g(a+J):= \sigma_g(a)+J$.

By Remark \ref{rem:from.Sira},
the exactness and 
residual proper outerness of 
$\sigma\colon G\to \mathrm{Aut}(A)$
allow
a natural identification
$$
(A\rtimes_{\sigma, \lambda}  G)/I  = 
(A/J)\rtimes _{\alpha  , \lambda}  G
\,.$$

Since $D\not\subseteq I $ 
implies $D_+ \not\subseteq I $, there exists
$d\in D_+$, $d\notin I $. 
The epimorphism 
$\pi_I\colon A\rtimes_{\sigma, \lambda}  
G\to (A\rtimes_{\sigma, \lambda}G)/I$ 
is equal to the quotient 
map $\pi^J$  from 
$A\rtimes_{\sigma, \lambda}  G$ 
onto 
$(A/J)\rtimes_{\alpha  , \lambda}  G \cong 
(A\rtimes_{\sigma, \lambda}  G)/ I $
(under natural identifications).
We denote the conditional expectation 
$E_\lambda \colon 
(A/J)\rtimes_{\alpha  , \lambda}  G \to A/J$
(temporarily) by $E$ and define 
$$b := \pi_I (d)\,,\,\quad   \text{and}  \,\quad 
\varepsilon : =\frac{1}{4} \,\|E(b)\|\, >0\,.$$
Lemma \ref{lem:SharpOlePed3}(i) 
gives an element
$x\in (A/J)_+$ such that
$$
\|x\|=1, \ \  \|xbx-xE(b)x\|<\varepsilon, \ \  \
\|xE(b)x\|>\, \|E(b)\|-\varepsilon \,= \,
\frac{3}{4}\, \| E(b)\|
\,.$$
By Remark \ref{rems:b-eps=d*ad}(i),  
there is a contraction 
$y\in (A/J)\rtimes _{\alpha , \lambda}  G$ 
such that 
$$ y^*xbxy\, =\,  
(xE(b)x-\varepsilon)_+\,  \in (A/J)_+\,.$$
Note that  $y^*xbxy\, \not=0\,$,  because
$$
\|(xE(b)x-\varepsilon)_+\|= 
\|xE(b)x\|-\frac{1}{4}\, \|E(b)\| 
> \frac{1}{2} \,\| E(b)\| = 2\varepsilon \,.
$$
Since 
$\pi_I | A= \pi^J |A$ and 
$(xE(b)x-\varepsilon)_+\in (A/J)_+$, 
there is
exists $c\in A_+$
such that 
$\pi^J(c)=(xE(b)x-\varepsilon)_+$.  
Since $\pi^J$ ($=\pi_I$) is surjective,
there exists a contraction 
$w\in A \rtimes _{\sigma, \lambda}  G$
with $\pi^{J}(w)=xy$. 
We obtain that 
$$c=w^*dw+ v $$
for some $v \in I$. 
The set $\Cf _c(G,J)$ is dense in $I$, 
because 
$I =J\rtimes_{\sigma, \lambda}  G$ and 
$G$ is discrete. 
This allows us to see, that  
$J I J$ is dense in $I$.
It follows that  
$ \{ e\in J_+\,; \,\, \| e \| <1 \} $
 is  an  approximate 
unit of $I $. 
In particular, there exists $e\in J_+$ with 
$\| v- ev \| < \varepsilon$.

Let $1$ denote the unit of 
$\widetilde{A}\rtimes_{\sigma, \lambda}  G$, 
then 
$A\rtimes_{\sigma, \lambda}  G$ is an ideal of 
$\widetilde{A}\rtimes_{\sigma, \lambda}  G$.
With 
$g:=(1-e)\in \widetilde{A}_+$, $\| g \| \leq 1$ 
we get 
$$
\|gw^*dwg  -  gcg\| = \|g v g\| \leq  \| v-ev \|< 
\varepsilon= \frac{1}{4}\, \| E(b) \|
\,.$$
Since  
$gzg=z+eze-(ze+ez)\,$ and 
$\pi_I (e)=\pi^J(e)=0$,
we have $\pi_I (gzg)=\pi_I (z)$ for all 
$z\in A\rtimes_{\sigma, \lambda}  G$.

By Remark \ref{rems:b-eps=d*ad}(i), 
there exists a contraction 
$h\in A\rtimes_{\sigma, \lambda}  G$ 
such that 
$$
h^* (gw^*dwg) h = 
(gcg-\varepsilon)_+   \in  A_+
\,.$$
With $z:=(d^{1/2}wgh)^*$ we have that 
$z^*z\in D$  and 
$z z^*= (gcg-\varepsilon)_+ =: f \in A_+\,$.
Finally, we see from $\pi_I (gcg)=\pi_I(c)$ that 
\begin{eqnarray*}
\|\pi_I (f)\|
&=&
\|\pi_I((gcg)-\varepsilon)_+)\|
=(\|\pi_I (gcg)\|-\varepsilon)_+=
(\|\pi_I (c)\|-\varepsilon)_+\\
&=&
(\|(xE(b)x-\varepsilon)_+\|-\varepsilon)_+ 
= \|xE(b)x\|
-\frac{1}{2}\| E(b) \|> \frac{1}{4} \| E(b) \| >0\,.
\end{eqnarray*} Hence, $f\not\in I$.
\qed
\section{Strongly purely infinite crossed products}
\label{sec:spi.crossed}
In this section we prove out main result 
Theorem \ref{thm:main}. We start with the 
definition of an $G$-separating action.
\begin{definition}
\label{def:G.separating}
Suppose that $(A,G,\sigma)$ is a 
\cst-dynamical system with discrete group $G$. 
The action of $G$ on $A$ is \emph{$G$-separating} 
if for every $a, b\in A_+\,$, $c\in A$, 
$\varepsilon>0$, there exist elements $s,t\in A$ 
and $g, h \in G$ such that
\begin{equation}
\label{InEq.G.separating} 
\| \,  s^*a\,s  - \sigma_g (a) \, \|  \, 
<\varepsilon\,,
\,\; \| \,   t^*b\,t  - \sigma_h (b)   \, \|  \, 
<\varepsilon \;\;
\textrm{and}\;\; \|  \, s^*c\,t\, \|  \, 
<\varepsilon\, .
\end{equation}
\end{definition}

\begin{rems}\label{G.separating}
(i)
Notice that Definition \ref{def:spi} 
and Remark \ref{rem:spi}
immediately implies that every action
$\sigma\colon G\to \mathrm{Aut}(A)$ is 
$G$-separating 
\emph{if $A$ itself is strongly purely 
infinite}: Take $h=g=e\in G$. 
If the contractions $s,t\in A$
satisfy the defining inequalities 
(\ref{InEq.spi2})
of strongly p.i.\  algebras $A$ then 
they also satisfy the
inequalities (\ref{InEq.G.separating}).
\smallskip

\noindent 
(ii)
$G$-separating 
actions on a locally compact  space 
$X$ are not necessarily minmal.
One can show that above mentioned
example of an  exact and non-minimal
action of the (small) Thompson group
$F$ on two parallel lines  
$\R \cup (i+\R) \subseteq \C$
is also $F$-separating.
\smallskip

\noindent 
(iii)
The existence of a $G$-separating action 
on $A$ imposes requirement on $A$ itself,
\eg 
the cases $a=b=c=p$ 
and $a=b=c=1$ 
with $\varepsilon=1/3\,$
in inequalities (\ref{InEq.G.separating})
show
that  \emph{$A$ can not contain
minimal  non-zero projections $p\in A$
and that
$1_A$ must be properly
infinite in $A$ if $A$ is unital.}
Therefore,  \cst-algebras, 
that are commutative \emph{and} unital,
can not have any
$G$-separating actions.
\smallskip

\noindent 
(iv)
Further variations of the 
concepts
that we introduce here are possible,
\eg one could start with
conditions that are weaker than 
conditions for
$G$-separating
actions.
Also  one  could require the existence
of $n\in \N$  such that 
for  $a,b\in A_+$
and $\varepsilon>0$ 
there is a solution 
$d_1,\ldots, d_n \in A$ and 
$g_1,\ldots , g_n \in G$ 
of the inequality 
(\ref{InEq.n.majorizing}) 
in
Definition \ref{def:n.majorizing.action}
of $n$-majorizing actions
whenever $b$
is in the smallest $G$-invariant 
closed ideal that contains $a$.
Or one could attempt to
replace
the filling family $\cF :=A_+$ by smaller
filling families 
$\cF \subseteq A_+$ and require more
elaborate local matrix diagonalization formulas
involving also $G$-translates, \cf 
Definition \ref{def:MatrixDiag}.
\end{rems}

Combing Theorem \ref{thm:A-filling} 
with
Theorem \ref{prop:local-spi}
we obtain the following result
\begin{thm}\label{cor:A-filling}
Let $G$ be a discrete 
group acting by 
$\sigma\colon G\to \mathrm{Aut}(A)$ 
on a separable or commutative \cst-algebra $A$. 
Suppose that the action 
is residually properly outer 
(\cf Def.\ \ref{def:properly.outer})
and exact
(\cf Def.\ \ref{def:exact}).
Let $\cF\subseteq  A_+$ be a filling family for $A$.
Then the following are equivalent:
\begin{itemize}
\item[(i)] The crossed product 
$A\rtimes_{\sigma, \lambda}  G$ 
is strongly purely infinite.
\item[(ii)] The family $\cF$ has the 
diagonalization property in 
$A\rtimes_{\sigma, \lambda}  G$.
\end{itemize}
\end{thm}
\pf Let $B:=A\rtimes_{\sigma, \lambda}  G$.
The assumptions ensure that 
$A_+$ is a filling family
for $B$ by Theorem \ref{thm:A-filling}.
Since $\cF$ is filling for $A$, 
$\cF$ is also filling for $B$
by Lemma 
\ref{lem:F.fill.A.A+.fill.B.Then.F.fill.B}.

(i)$\Rightarrow$(ii): If $B$  is s.p.i., 
then 
$B_+$ 
has the diagonalization property 
(see Definition \ref{def:mdiag.family})
in $B$, 
\cf  \cite[lem.~5.7]{KirRorOinf}.
This implies that our family 
$\cF\subseteq A_+\subseteq B_+$ 
has the diagonalization property
in $B$.

(ii)$\Rightarrow$(i): Since our family 
$\cF\subseteq  A_+$ is filling for 
$B$,
and since $\cF $ has the diagonalization
property in $B$,
we get that $B$ is s.p.i.\  by 
Theorem \ref{prop:local-spi}.
\qed

\begin{rem}\label{rem:e.is.dU}
Let $(A,G,\sigma)$ a \cst-dynamical system. 
\begin{itemize}
\item[(i)] For each 
$a_1,a_2 \in A_+$, $x,d_1,d_2\in A$, 
$g_0,g_1,g_2\in G$ 
and 
$s_1 := d_1U(g_1)$, 
$s_2:=\sigma_{g_0^{-1}}(d_2)U(g_0^{-1}g_2g_2)$,
$c:=xU(g_0)$,
$b_1:=a_1$, and $b_2:=\sigma_{g_0}(a_2)$
the following 
equalities hold:
$$\|  s_j^*a_js_j - a_j\| = 
\| d_j^*b_j d_j - \sigma_{g_j}(b_j)\|\,
\quad
\textrm{and}\,
\quad
\| s_1^*cs_2\|=\| d_1^*x d_2\|
\,.$$	
\item[(ii)] With $g_0=e$ in (i) the equalities 
reduce to:
$$\|  s_j^*a_js_j - a_j\| = 
\| d_j^*a_j d_j - \sigma_{g_j}(a_j)\|\,
\quad
\textrm{and}\,
\quad
\| s_1^*cs_2\|=\| d_1^*c d_2\|
\,.$$
\end{itemize}
\end{rem}

\begin{prop}\label{prop:infinite}
Suppose that $(A,G,\sigma)$ is a 
\cst-dynamical system with discrete $G$. 
The following properties \emph{(i)--(ii)} 
are equivalent:

\begin{itemize}
\item[(i)] The action of $G$ on $A$ is 
$G$-separating
in the sense of Definition 
\ref{def:G.separating}.
\item[(ii)] For every $a_1,a_2\in A_+$, 
$c\in A\rtimes_{\sigma, \lambda}  G$ 
and $\varepsilon>0$ there exist 
$d_1,d_2\in A$ and
$g_1,g_2\in G$
such that the elements 
$s_j=d_jU( g_j)$ of $\Cf _c(G,A)$ satisfy, for 
$j=1,2$,
\begin{eqnarray}\label{eqn1.a}
\| \,  s_j^*a_js_j-a_j\,\|  \, <\varepsilon\,,\;
\;\;\textrm{and}\;\;  \| \, s_1^*cs_2\,\|  \, <
\varepsilon\, .
\end{eqnarray}
\end{itemize}
\end{prop}
\pf
(ii)$\Rightarrow$(i): 
If we take 
$c\in A$,  $a_1:=a$ and $a_2:=b$
for $a,b\in A_+$  
and $\varepsilon>0$, then (ii) implies, 
using Remark \ref{rem:e.is.dU}, that there exist 
$d_1,d_2\in A$ and
$g_1,g_2\in G$
such that 
$\|d_j^*a_j d_j - \sigma_{g_j}(a_j)\|<\varepsilon$ 
and $\|d_1^*cd_2\|<\varepsilon$, so
the inequalities
(\ref{InEq.G.separating}) of  
Definition \ref{def:G.separating}
are satisfied 
with $d_1,d_2,g_1,g_2$ in 
place $s,t,g,h$. 

(i)$\Rightarrow$(ii): 
Define 
$\Co:=\,\{ dU(g) \,;\,\, d\in A\,,\,\, g\in G\,\}$
and 
$\cS := \Co \,$. 
Select any $\varepsilon_0>0$.
Clearly, the closed linear span of 
$\Co$ is equal to
$A\rtimes_{\sigma, \lambda} G$. 
If we can show that
$\cF:=A_+\,$, $\Co$ and
$\cS$ satisfy the assumptions (i)-(iii) of 
Lemma \ref{lem:cS.control} -- with 
$A\rtimes_{\sigma, \lambda}  G$ in place of 
$A$
--, then it
follows from Lemma \ref{lem:cS.control} 
that for every $a_1,a_2\in A_+$, 
$c\in A\rtimes_{\sigma, \lambda}  G$ and 
$\varepsilon>0$ there exist $d_1,d_2\in A$ 
and
$g_1,g_2\in G$
such that
$s_j=d_jU( g_j)\in \cS$ fulfil \eqref{eqn1.a},
which in turn gives (ii).

It is evident that 
our $\Co$ and $\cS$
satisfy properties  (ii) and (iii)
in Lemma \ref{lem:cS.control}. 
Since $A_+$ is 
closed under $\varepsilon$-cut-downs,
property (i) becomes automatic if 
each pair $(a_1,a_2)$, with 
$a_1,a_2\in A_+$, 
has the matrix diagonalization property 
of Definition \ref{def:MatrixDiag}
with 
respect to $\cS$ and $\Co$:

If $a_1,a_2\in A_+$, 
$c=xU(g_0)\in \Co$ 
with $x\in A$, $g_0\in G$, and $\varepsilon>0$
are given, then we define
$b_1:=a_1$, $b_2:=\sigma_{g_0}(a_2)$. 
(If we instead of $\varepsilon$ are 
given positive 
$\varepsilon_1$, $\varepsilon_2$ and $\tau$, 
set $\varepsilon:=
\min(\varepsilon_1,\varepsilon_2,\tau)$.) 
Since the action $\sigma$ is $G$-separating, 
we can find $d_1,d_2\in A$ and $g_1,g_2\in G$ with 
$\|d_1^*b_1d_1-\sigma_{g_1}(b_1)\|$,  
$\|d_2^*b_2d_2-\sigma_{g_2}(b_2)\|$
and $\| d_1^*xd_2 \|$ all strictly below 
$\varepsilon$.
Remark \ref{rem:e.is.dU} 
provides elements $s_j\in \Co$
satisfying \eqref{eqn1.a}.
Thus
$(a_1,a_2)$ has the matrix 
diagonalization property with 
respect to $\cS$ and $\Co$.
\qed

\begin{thm}\label{thm:F-filling}
Let 
$(A,G,\sigma)$  a \cst-dynamical system,
with discrete  $G$. 
Suppose that $A_+ $ is a filling 
family for 
$A \rtimes_{\sigma, \lambda}   G$
and that the action of $G$ on $A$ is 
$G$-separating. 
Then 
$A\rtimes_{\sigma,\lambda} G$ 
is strongly purely infinite.
\end{thm}

\pf
By Theorem \ref{prop:local-spi} it
remains to show that $A_+$ has the
diagonalization property 
in $A \rtimes_{\sigma, \lambda}   G$.
Since $A_+$ is closed under 
$\varepsilon$-cut-downs 
Lemma \ref{lem:flexible.diag.imply.general}
applies, and therefore it is enough to show
that each pair
$(a_1,a_2)$ with $a_1,a_2\in A_+$,
has the matrix diagonalization property 
in $A \rtimes_{\sigma, \lambda}   G$.
But this follows from the $G$-separation 
property of the action $\sigma$ by 
Proposition  \ref{prop:infinite}.
\qed.

\textbf{Proof  of Theorem \ref{thm:main}}:\, 
By Theorems 
\ref{thm:A-filling} and \ref{thm:F-filling}
the assumptions imply that
$A\rtimes_{\sigma, \lambda}  G$ 
is strongly purely infinite.\qed

\begin{rem}\label{rem:general.}
Suppose that $(A,G,\sigma)$ is a 
\cst-dynamical system and 
that $G$ is discrete. Then a
family  
$\cF \subseteq  
A_+  \subseteq  A\rtimes_{\sigma, \lambda}  G\,$ 
which is \emph{invariant under 
$\varepsilon$-cut-downs} 
has the 
diagonalization properly in 
$A\rtimes_{\sigma, \lambda}  G\,$, 
if and only if, for every $a_1,a_2\in \cF$, 
$c \in \Cf _c(G,A)$
and $>0\,$, there exist 
$s_1, s_2\in \Cf _c(G,A)$ such that, for 
$j=1,2$,
\begin{equation*}
\| \,  s_j^*a_js_j-a_j\,\|  \, <\varepsilon\,
\quad
\textrm{and}\,
\quad
\| \, s_1^*cs_2\,\|  \, <\varepsilon\,.
\end{equation*}
This follows from
Lemma \ref{lem:flexible.diag.imply.general},
Lemma \ref{lem:cS.control} and 
the fact that $\Cf _c(G,A)$ is dense 
in $A\rtimes_{\sigma, \lambda}  G\,$.
\end{rem}

\begin{rem}\label{rem:trivial.action.of.exact.G}
Notice that for an 
exact locally compact group $G$ 
the reduced group 
\cst-algebra $C^*_\lambda(G)$ 
is an exact \cst-algebra 
(\cf
\cite[p.~171]{KirWa.Exact.Groups.Cont.Bund}).
By Theorem \ref{thm:A.ot.B.spi.if.A.spi.B.exact}, 
the minimal \cst-tensor product 
$A\otimes^{\min}B$ of a s.p.i.\ \cst-algebra 
$A$ with an exact \cst-algebra
$B$ is again s.p.i. 
Hence, if $G$ is an exact locally compact group,
$\sigma(g):=\id_A$ is the trivial action
on a s.p.i.\ \cst-algebra $A$ then 
$A\rtimes_{\sigma,\lambda} G \cong  
A\otimes^{\min} C^*_\lambda (G)$
is s.p.i.

This shows that there is 
room for 
refinements of our
sufficient conditions on the actions that imply 
strong pure infiniteness of the 
reduced crossed products:
Here the action 
$\sigma$ is even not 
element-wise properly outer,
but satisfies the 
$G$-separation property and is an
exact action by \cst-exactness of 
$C^*_\lambda (G)$.
\end{rem}

\section{The case of commutative \cst-algebras}
\label{sec:abelian.case}
The case of  
$G$-actions on commutative \cst-algebras
allows some topological interpretation. 
The next lemma  
has inspired our choice of 
the name \emph{$G$-separating} 
in Definition 
\ref{def:G.separating}.
Notice that we do not require 
$\alpha_{g_j}(U_j)\subseteq U_j$
in its part (ii).
\begin{lem}\label{lem:G-separating}
Suppose that 
$(A,G,\sigma)$ is a \cst-dynamical system, 
that 
$A\cong \Cf _0(X)$ is commutative, 
and that 
the action $\sigma$ of $G$ on  
$\Cf _0(X)$  is induced 
by the action $\alpha$ of $G$ on 
$X\cong \widehat{A}$
(\ie 
$\sigma_g(f):=f\circ \alpha_{g}^{-1}$ 
for $f\in A, g\in G$). 
Then the following properties are equivalent:
\begin{itemize}
\item[(i)] The action of $G$ on $A$ is 
$G$-separating, \ie  for 
every $a, b\in A_+\,$, $c\in A$, $\varepsilon>0$, 
there exist elements 
$d_1,d_2\in A$ and $g_1,g_2\in G$ such that
$$
\| \,  d_1^*ad_1-\sigma_{g_1}(a)\,\|  \, <
\varepsilon\,,\;
\| \,  d_2^*bd_2-\sigma_{g_2}(b)\,\|  \, <
\varepsilon
\;\;\textrm{and}\;\; \| \, d_1^*cd_2\,\|  \, <
\varepsilon\, .
$$
\item[(ii)] For every open 
$U_1,U_2\subseteq X$  and 
compact $K_1,K_2\subseteq X$ with 
$K_1\subseteq U_1$, $K_2\subseteq U_2$,
there exist $g_1,g_2\in G$ such that 
$$
\alpha_{g_1}(K_1)\subseteq U_1, \ \ 
\alpha_{g_2}(K_2)\subseteq U_2,
 \ \ \alpha_{g_1}(K_1)\cap \alpha_{g_2}(K_2)
 =\emptyset
 \,.$$
\end{itemize}
\end{lem}
\pf 
(ii)$\,\Rightarrow$(i):
Let $a,b\in A_+$, $c\in A$ and 
$\varepsilon >0$. 
We use  assumption (ii) on
$$
U_1:=a^{-1}(\varepsilon/4,\infty)=
\{x\in X \,;\,\, a(x)>\varepsilon/4\} \, ,\, \,\, \,
U_2:=\{x\in X \,;\,\, b(x)>\varepsilon/4\}\, ,
$$
$$
K_1:=\{x\in X \,;\,\, a(x)\geq \varepsilon/2\}\,, 
\, \,\, \, 
K_2:=\{x\in X \,;\,\, b(x)\geq \varepsilon/2\}\,,$$
and find $g_1, g_2\in G$
such that 
$$
\alpha_{g_1}(K_1)\subseteq U_1\,, \ \ 
\alpha_{g_2}(K_2)\subseteq U_2\,,
 \ \ \alpha_{g_1}(K_1)\cap 
 \alpha_{g_2}(K_2)=\emptyset
 \,.$$
Since $a,b\in \Cf_0(X)_+$, 
we have that 
$\overline{U_1} \subseteq 
a^{-1}[\varepsilon/4,\infty)$
and 
$\overline{U_2} \subseteq 
b^{-1}[\varepsilon/4,\infty)$ 
are  compact subsets of $X$. 

Since the compact sets 
$\alpha_{g_1}(K_1)$ and 
$\alpha_{g_2}(K_2)$ 
are disjoint, 
applications of Tietze extension theorem
gives
elements $e_1,e_2\in A_+$
with $\| e_j\| \leq  2/\sqrt{\varepsilon}$ 
and a contraction $f=f^*\in A$ such that
$$
e_1 | {\overline{U_1}}=
a^{-1/2} | {\overline{U_1}}\,, \,\, \, \,
e_2 | {\overline{U_2}}=
b^{-1/2} | {\overline{U_2}}\,, \, \, \, \,
f | {\alpha_{g_1}(K_1)}= -1\,, 
\, \, \, \,
f | {\alpha_{g_2}(K_2)}=1
\,.$$
Let $f_+,f_-\in A_+$ be
the canonical decomposition 
$f=f_+-f_-$ with 
$f_+f_-=0$.
We define
 $$
 d_1:=\, 
 e_1(\sigma_{g_1}(a)-\varepsilon/2)_+^{1/2} f_-\, 
 \,\,\,\,\, \text{and} \,\,\,\,\,\,
 d_2:=\, 
 e_2 (\sigma_{g_2}(b)-\varepsilon/2)_+^{1/2} f_+
 \,.$$
Then $d_1^*cd_2=0$ because $f_+f_-=0$.

Since 
$(\sigma_{g_1}(a)-\varepsilon/2)_+(x)\not=0$ 
implies 
$a(\alpha_{g_1^{-1}}(x))>\varepsilon/2$, 
we get 
$\alpha_{g_1^{-1}}(x)\in 
K_1$, and
$x\in \alpha_{g_1}(K_1)\subseteq 
U_1\subseteq \overline{U_1}$. 
It implies that
$f_-(x)=1$ and $e_1(x)=a^{-1/2}(x)$.
We obtain that
$$
d_1^*ad_1=
e_1^2a(\sigma_{g_1}(a)-\varepsilon/2)_+(f_-)^2=
(\sigma_{g_1}(a)-\varepsilon/2)_+
\,.$$
In a similar way we see that 
$d_2^*bd_2=(\sigma_{g_2}(b)-\varepsilon/2)_+$. 

(i)$\,\Rightarrow$(ii): 
Let 
$U_1,U_2\subseteq X$ be open and 
$K_1,K_2\subseteq X$ compact subsets such that 
$K_1\subseteq U_1$ and $K_2\subseteq U_2$. 
We can assume that the intersection 
$K_1\cap K_2$ is non-empty.
There exists an open set 
$W$ with a compact closure 
$\overline{W}$ such that
$$K_1\cup K_2\, \subseteq \,
W\, \subseteq 
\overline{W}\, \subseteq \,U_1\cup U_2\,.$$
By the Tietze extension theorem, 
there are contractions 
$a,b,c\in A_+$ such that
$$
a | {K_1}\, =1\,, \, \,\,  b | {K_2}\, =1\,, 
\, \, \,  c | {\overline{W}}\, =1
\,, $$
$$
\mathrm{supp}(a)\subseteq 
U_1\cap W, \ \ \mathrm{supp}(b)\subseteq 
U_2\cap W, \ \ \mathrm{supp}(c)\subseteq 
U_1\cup U_2
\,.$$
We apply assumption (i) to  $a$, $b$, $c$ and 
$\varepsilon:=1/4$, and 
get elements $d_1,d_2\in A$ and 
$g_1,g_2\in G$ such that
$$\|d_1^*ad_1-\sigma_{g_1}(a)\|< 1/4\,, \, \,\, 
\|d_2^*bd_2-\sigma_{g_2}(b)\|< 1/4\,, \,\,\,  
\|d_1^*cd_2\| < 1/4\,.$$
If  $x\in \overline{W}$, then  $c(x)=1$ and
$|d_1(x)||d_2(x)|
\leq\|d_1^*cd_2\|<1/4$. Thus,
$$
V_1:=\{x\in \overline{W}\,;\,\,\, 
|d_1(x)|\geq 1/2\}\,, 
\, \,\, 
V_2:=\{x\in \overline{W} \,;\,\, 
|d_2(x)|\geq 1/2\}
$$
are disjoint sets. 
If $x\in \alpha_{g_1}(K_1)$, then 
$\alpha_{g_1^{-1}}(x)\in K_1$
and 
$\sigma_{g_1}(a)(x)=a(\alpha_{g_1^{-1}}(x))=1$.
It follows $|d_1(x)|^2 a(x) \ge 3/4$ 
and $|d_1(x)|>1/2$
(the latter because $1\ge a(x) >0$). 
Thus, $x\in U_1\cap W$
and  $x\in V_1$. It shows 
$\alpha_{g_1}(K_1)\subseteq U_1\cap V_1$.
In a similar way we get  
$\alpha_{g_2}(K_2)\subseteq U_2\cap V_2$. 
It implies 
 $\alpha_{g_1}(K_1)\subseteq U_1\,$, 
$\,\alpha_{g_2}(K_2)\subseteq U_2$ 
and that 
$\alpha_{g_1}(K_1)\cap \alpha_{g_2}(K_2)
=\emptyset\,$.\qed

\medskip

The following condition (i) in 
Corollary \ref{cor:Abelian}
is satisfied if the action $\alpha$ 
has the residual version
of the topological freeness in sense of 
\cite[def.~1]{ArchSpiel},
see \eg  the essential freeness defined in 
\cite[def.~1.17]{Siera2010}
(inspired by
 \cite[def.~4.8]{Ren:fixed}) when 
 $G$ is countable.

\begin{cor}\label{cor:Abelian} 
Let $G$ be a discrete group,
$\alpha\colon G \to \mathrm{Homeo}(X)$ 
an action of $G$
on a locally compact Hausdorff space $X$. 
Suppose that 
\begin{itemize}
\item[(i)]  for every closed 
$G$-invariant subset $Y$ of $X$ and
every  $e\not= g\in G$ 
the set $\{y \in Y\colon \alpha_g(y)=y\}$ 
has empty interior,
\item[(ii)]   the action 
$\sigma\colon G \to \mathrm{Aut}(\Cf_0(X))$,
given by $\sigma_g(f):=f\circ (\alpha_{g})^{-1}$,
is exact on the \cst-algebra $\Cf_0(X)$,
and
\item[(iii)] 
the action $\alpha$
is $G$-separating, 
\ie by Lemma \ref{lem:G-separating},
for every $U_1,U_2\subseteq X$ open and 
$K_1,K_2\subseteq X$ 
compact such that $K_1\subseteq U_1$, 
$K_2\subseteq U_2$, 
there exist $g_1,g_2\in G$ such that 
$$
\alpha_{g_1}(K_1)\subseteq 
U_1, \ \ \alpha_{g_2}(K_2)\subseteq U_2,
 \ \ \alpha_{g_1}(K_1)\cap \alpha_{g_2}(K_2)=
 \emptyset
 \,.$$
\end{itemize}
Then 
$\Cf _0(X)\rtimes_{\sigma, \lambda}  G$ 
is a strongly purely infinite \cst-algebra.
\end{cor}
\pf
Let $A:=\mathrm{C}_0(X)$.
It is easy to see that property (i) implies 
that the action on any quotient $A/I$ by a 
$G$-invariant closed ideal $I \neq A$ is 
element-wise properly outer. 
Now Theorem \ref{thm:main} applies to 
$A\rtimes _{\sigma, \lambda}  G$.
\qed

The following remark shows that
in case of commutative $A$ and discrete
\emph{amenable} $G$ several of the
previously  considered properties are
equivalent.
\begin{rem}\label{rem:converse.direction}
If $A$ is commutative and $G$ is a 
discrete \emph{amenable} group 
that acts
on $A$ by $\sigma$, 
then the following properties are
equivalent:
\begin{itemize}
\item[(i)] $A$ separates the ideals of 
$A\rtimes _{\sigma} G$, \ie
$I\mapsto A\cap I$ is an injective map from 
$\mathcal{I }(A\rtimes _{\sigma} G)$ into 
$\mathcal{I}(A)$
(see \cite[def.~1.9]{Siera2010}).
\item[(ii)] The action 
$\sigma\colon G\to \mathrm{Aut} (A)$ 
is residually
properly outer (Definition  
\ref{def:properly.outer}).
\item[(iii)] 
The family
$\cF:=A_+$ is filling for 
$A\rtimes _{\sigma} G$ 
(Definition \ref{def:filling.mdiag.family}).
\end{itemize}
%
\pf
(i)$\Rightarrow$(ii): 
By \cite[thm.~4.1]{KawamuraTomiyama} 
(in case of unital $A$, and 
\cite[thm.~2]{ArchSpiel} for the general
-- non-unital --  case) the separation property
implies that the adjoint action of $G$ on the
Gelfand spectrum of $A$ is topologically free, 
which is equivalent to element-wise proper
outerness by \cite[prop.~1]{ArchSpiel}.

This applies also to the quotients 
$(A/J)\rtimes_{[\sigma]_J, \lambda} G$,
because the property (i) 
passes to quotients by amenability
of $G$.
See also \cite[thm.~1.13]{Siera2010}.

(ii)$\Rightarrow$(iii):
Since amenable  $G$ are exact, 
the residual proper outerness of
the action implies that 
$\cF:=A_+$ is filling for 
$A\rtimes _{\sigma,\lambda} G$
by Theorem \ref{thm:A-filling}.

(iii)$\Rightarrow$(i):
By Remark \ref{rems:filing.familiy}, 
the subalgebra $A$ separates the ideals of 
$B:= A\rtimes _{\sigma,\lambda} G$
if $\cF:= A_+$ is filling for $B$.\qed

Asking $G$ to be amenable can be weakened 
to exactness of $\sigma$ and 
$A\rtimes _{\sigma,\lambda} G\cong 
A\rtimes _{\sigma} G$. One might also expect 
nuclearity of $A\rtimes _{\sigma,\lambda} G$
would suffice in place of amenability of $G$ 
(this is know at least in the unital case).
\end{rem}

\section{Strong boundary actions 
versus {\it G}-separating actions}
\label{sec:general.strong.boundary.actions}
In this section we prove our 
Theorem \ref{thm:n.majorizing.n.covering}. 
We state 
with the definition of $n$-majorizing 
and $n$-covering actions.
\begin{definition}
\label{def:n.majorizing.action} 
Let $n\in\N$ and  $A$ a non-zero 
\cst-algebra, that is not isomorphic to a 
subalgebra of $M_{n+1}(\C)$ 
\emph{if $A$ is unital}. An action 
$\sigma\colon G\to \mathrm{Aut}(A)$ will be 
called an \emph{$n$-majorizing} action of 
$G$ on $A$, if, for every non-zero $a\in A_+$, 
every \emph{non-invertible} $\,b\in A_+\,$ 
and every $\varepsilon >0$, there exist 
$d_1,\ldots, d_n \in A$ and 
$g_1,\ldots , g_n \in G$ such that
\begin{equation}\label{InEq.n.majorizing}
\| \,\sum_{j=1}^n \, d_j^*\, \sigma_{g_j} 
(a) \,d_j\, - \, b\, \| \, <\, \varepsilon\,.
\end{equation}
\end{definition}
\begin{definition}\label{def:n.covering.action} 
Let $n\in \N$, $n\ge 2$ 
(\,\footnote{\, If $n=1$ then property 
(\ref{InEq.n.covering}) holds if and only if 
$A$ is a unital simple purely infinite 
\cst-algebra.}\,). 
Suppose that $A$ is a unital \cst-algebra, that 
is not isomorphic to a *-subalgebra of $M_{n}(\C)$. 
An action $\sigma$  of $G$ on $A$ is an 
\emph{$n$-covering action}  if, for every non-zero 
$a\in A^+$, and every $\varepsilon  >0$, there 
exist $d_1,\ldots, d_n \in A$ and 
$g_1,g_2,\dots,g_n\in G$ and such that 
\begin{equation}
\label{InEq.n.covering} 
\, \| \,  \sum_{j=1}^n \, d_j^*\, \sigma_{g_j} 
(a) \,d_j\, - \, 1\, \| \, <\, \varepsilon\,.
\end{equation}
\end{definition}
The following lemma denies the existence of 
non-zero ``socles'' in \cst-algebras $A$ 
that admit $n$-majorizing or $n$-covering actions
considered in Definitions
\ref{def:n.majorizing.action} 
and \ref{def:n.covering.action}. 
\begin{lem}
\label{lem:deny.socle.for.major.cover}
Let  $(A,G,\sigma)$
\cst-dynamical system.  

Consider the following properties 
\emph{($\alpha$)} or \emph{($\beta$)} of 
$(A,G,\sigma)$ depending on $n\in \N$: 
\begin{itemize}
\item[($\alpha$)]
There is $n\ge 1$ such that,
for each non-zero $a\in A_+$,
\emph{non-invertible} 
$b\in A_+$ and $\varepsilon>0$, 
there exist
$d_1,\ldots, d_n \in A$ and 
$g_1,\ldots , g_n \in G$ 
that satisfy the inequality 
\emph{(\ref{InEq.n.majorizing})}
in Definition 
\ref{def:n.majorizing.action}.
\item[($\beta$)]  $A$ is unital, and
there is $n\ge 2$ such that,
 for each non-zero 
 $a\in A_+$ and $\varepsilon>0$,
there exist 
$d_1,\ldots, d_n \in A$ and 
$g_1,\ldots , g_n \in G$ 
that satisfy the inequality 
\emph{(\ref{InEq.n.covering})}
in Definition \ref{def:n.covering.action}.
\end{itemize}
If $A$ is unital and $(A,G,\sigma)$ satisfies 
($\alpha$)
then it satisfies ($\beta$) with
$n$  replaced by $n+1$, and if $(A,G,\sigma)$
satisfies ($\beta$) then it satisfies ($\alpha$) 
-- with same $n\in \N$.
If $(A,G,\sigma)$ satisfies 
($\alpha$) or ($\beta$), then
the algebra $A$ is $G$-simple, 
\ie
$A$ and $\{ 0 \}$ are the only 
$G$-invariant closed
ideals of $A$.

If $A$ contains
a projection $p\not=0$ with 
$pAp=\mathbb{C}\cdot p$,
then $A$ is a \cst-subalgebra of $M_{n+1}$ 
(respectively of $M_n$)
if  $(A,G,\sigma)$ has property ($\alpha$)
(respectively has property ($\beta$)).

The shift action 
$\sigma$ of the cyclic group 
$\Z_{n+1}$ on $A:=\Cf (\Z_{n+1})$
satisfies \emph{($\alpha$)} for 
$n\in \N$ 
and is element-wise properly outer.
\end{lem}
\pf
If $A$ is unital and 
$\sigma$ has property ($\alpha$), 
then let
$b:= 1-(\|a\|^{-1} a)^3$, and find 
$d_1,\ldots,d_n$ and 
$g_1,\ldots, g_n$ that
satisfy the inequality 
(\ref{InEq.n.majorizing}).
If we let  $g_{n+1}:=e$  and  
$d_{n+1}:= \|a\|^{-3/2} a$,
then $a$ and 
$g_1,\dots, g_n,g_{n+1}$ satisfy
(\ref{InEq.n.covering}).
It shows that actions on unital $A$
that satisfy property ($\alpha$)
are also actions that satisfy 
($\beta$) with $n+1$
in place of $n$.
If $(A,G,\sigma)$ satisfies 
($\beta$) and non-zero
elements $a,b\in A_+$ are given with 
$\| b\|=1$, then 
$d_1b^{1/2},\ldots,d_nb^{1/2}$
and $g_1,\ldots, g_n$ is a solution
of the inequality (\ref{InEq.n.majorizing})
in Definition \ref{def:n.majorizing.action}  
if $d_1,\ldots,d_n$ and $g_1,\ldots, g_n$
satisfy the
inequality (\ref{InEq.n.covering})
of Definition \ref{def:n.covering.action}.

The properties ($\alpha$) and ($\beta$)
imply  that $\{ 0\}$ and $A$ are the only 
$G$-invariant closed ideals
of $A\,$:\,  
If $A$ is \emph{non-unital} in case 
($\alpha$),
then the approximation, 
as expressed by the inequalities 
(\ref{InEq.n.majorizing}),  
shows that for each non-zero
$a\in A_+$  the smallest closed 
$G$-invariant ideal of
$A$ containing $a$ contains each 
$b\in A_+$.
If $A$ is unital and the actions has property 
($\beta$)
then each $G$-invariant closed ideal of $A$
contains $1$.
If $A$ is unital and the
C*-dynamical system 
$(A,G,\sigma)$ satisfy
property ($\alpha$) 
then it satisfies
property ($\beta$) with $n$ 
replaced by 
$n+1$.
Thus, again, $A$ and 
$\{0\}$ are the only closed
$G$-invariant ideals.

\smallskip

From now on, we suppose that
there exists a projection 
$0\not=p\in A_+$  with 
$pAp=\mathbb{C}p$. We call
those projections ``minimal'', 
even if minimal non-zero
projections 
of a \cst-algebra $A$
do not have the property 
$pAp=\mathbb{C}p$ in general,
e.g.\ the unit of the
Jiang-Su algebra $\mathcal{Z}$
is a minimal projection.
We show that this assumption
of the existence of such $p\in A$, 
together with the assumption that
$\sigma$ satisfies 
($\alpha$), implies that $A$
is unital. Thus
$A$ satisfies ($\beta$)
\emph{with $\,n+1\,$ in place of $\,n\,$}.
Then we derive that property ($\beta$)
and the existence of such $p\in A$
imply that $A$ is a \cst-subalgebra of
$M_n$.

It is obvious
that the ideal $\textrm{socle}(A)$ 
generated by those 
``minimal''  projections 
is (universally)  invariant under all  
automorphisms of $A$, \ie  
$\alpha(\mathrm{socle}(A))= 
\mathrm{socle}(A)$ 
for all 
$\alpha\in \mathrm{Aut}(A)$.
This happens also for the closure 
$J$ of $\mathrm{socle}(A)$.
Thus, $J$ must be $G$-invariant.
It follows that $J=A$ using 
$\mathrm{socle}(A)\neq \emptyset$.

It is not difficult to see, 
that $J$ is isomorphic to the 
$c_0$-direct sum 
of a family of algebras 
$\cK(\Hi _\tau)$ of 
compact operators on
suitable Hilbert spaces $\Hi_\tau$, 
and that $p$ is a rank-one projection 
on  some $\Hi_{\tau_0}$.  
Let $\Hi$ denote the 
Hilbert space sum of the Hilbert
spaces $\Hi_\tau$.
Then  $A$ becomes isomorphic to a 
non-degenerate
\cst-subalgebra of the algebra of compact
operators $\cK(\Hi)$ on $\Hi$, 
in a way that each minimal non-zero 
projection $p\in A$ 
becomes a rank one projection
on $\Hi$. This happens also for all 
$\sigma_g(p)$.
Recall that every projection in 
$A\subseteq  \cK(\Hi)$ has finite rank in 
$\cL (\Hi)$. 
Since  $A$ is a \cst-subalgebra
of  $\cK(\Hi)$, $A$
is in particular
an AF-algebra, and 
-- therefore -- contains an approximate 
unit $(q_\tau)$  consisting of 
an upward directed net of projections 
of finite rank in $\Hi$.

We show that $A$ must be unital if 
$(A,G,\sigma)$ satisfies ($\alpha$) in 
addition: Suppose that $A$ is not unital, 
then none of the projections 
$(q_\tau)$ are invertible in $A$.
Therefore, we can take
$b:=q_\tau$, $a:= p$ and
$\varepsilon=1/2$ in 
(\ref{InEq.n.majorizing}).
It follows that each
$q_\tau$ has linear rank $\leq n$.
This implies that the approximate unit
$(q_\tau)$ 
must be constant
$q_\tau=e$ for all $\tau\ge \tau_1$
with suitable $\tau_1$.
Then $e\in A$ is necessarily
the unit element of $A$, in contradiction
to our assumption that $A$ 
is \emph{not} unital.

It follows that $A$ must be unital,
and -- as above
observed -- the action 
$\sigma$ satisfies property
($\beta$)  with $n$  replaced by $n+1$.

If $A$ is unital and $(A,G,\sigma)$ 
satisfies property $(\beta)$,
then we  take  $a:=p$ and 
$\varepsilon:= 1/2$
in inequalities (\ref{InEq.n.covering}).
It shows that the rank of 
$1_A$ in its representation is $\leq n$. 
Thus $A$ is a \cst-subalgebra of $M_n$
in case ($\beta$).

\smallskip

The crossed product   
$\Cf(\Z_{n+1})\rtimes \Z_{n+1} $  is
naturally isomorphic to $M_{n+1}$.  
Hence, by Remark
\ref{rem:converse.direction}, the action 
of $\Z_{n+1}$ on $\Cf(\Z_{n+1})$ is 
element-wise properly outer. 
If $a\in A_+:=\Cf(\Z_{n+1})_+$ 
is non-zero and 
$b\in A_+$ is not invertible, 
then there are non-zero minimal projections
$p,q\in A_+$ and 
$\delta>0$ such that
$a\ge \delta p$ and 
$b\leq  \| b \| \cdot (1-q)$. 
Select
$g_1,\ldots, g_n \in \Z_{n+1}$ with 
$\sum_{j=1}^n \sigma_{g_j}(p)=1-q$.
It implies that
$\delta^{-1} \sum_{j=1}^n \sigma_{g_j}(a) 
\ge 1-q$.
Thus, there exists 
$T\in (1-q)A_+(1-q)$  with
$T (\sum_{j=1}^n \sigma_{g_j}(a))T=1-q$. 
Then
$a$, $b$,  $d_j:= T b^{1/2}$, 
$j=1,\ldots,n$ and $g_1,\ldots, g_n$
satisfy the inequality 
(\ref{InEq.n.majorizing})
for each $\varepsilon>0$.
\qed
\begin{rem}\label{rem:simple.spi}
Let $B$ be a non-zero simple
\cst-algebra. In preperation for 
the proof of Theorem 
\ref{thm:n.majorizing.n.covering}
we display here a 
number of properties equivalent
to strong pure infiniteness of $B$.
\begin{itemize}
\item[(i)] $B$ is strongly purely infinite.
\item[(ii)] Each non-zero element of 
$B_+$ is properly infinite 
in sense of \cite{KirRor1}.
\item[(iii, $n$)] There exists $n\in \N$ such 
that, for each non-zero elements $a,b \in B_+$
and $\varepsilon>0$, there exists 
$n$ elements
$d_1,\ldots, d_n\in B$ with 
\begin{equation}\label{InEq.pi.n}
\| d_1^*a d_1+\cdots + d_n^*ad_n -b \| 
<\varepsilon\,,
\end{equation}
and $B$ is not isomorphic to 
$M_k$ for any $k\leq n$.
\item[(iv)] 
$B$ is locally purely infinite in sense of 
\cite[def.~1.3]{BlanKir2}, \ie each non-zero 
hereditary
\cst-subalgebra of 
$B$ contains a non-zero stable
\cst-subalgebra.
\item[(v)] $B$ is purely infinite 
in the sense of J.~Cuntz
\cite[p.~186]{Cuntz:K-Th.in.Annals}, 
\ie each non-zero
hereditary \cst-subalgebra contains 
an infinite projection.
\end{itemize}
\end{rem}

\pf Property (ii) implies (i) by 
\cite[thm.~5.8]{BlanKir2} and (i) 
implies (ii)  by \cite[prop.~5.4]{KirRorOinf}.
Property (iii,$n=1$) is equivalent to (ii) 
by \cite[thm.~4.16]{KirRor1}.
The properties  (iii,$n=1$),  
 (iv) and (v) are equivalent
by  \cite[prop.~3.1]{BlanKir2}.

(iii)$\Rightarrow$(ii): 
Suppose that $B$ is elementary, \ie
that  $B\cong \cK (\Hi)$ for some
Hilbert space $\Hi$.
By (iii), applied on some rank-one 
projection $a:=p\in B_+$, $b\in B_+$, and 
$\varepsilon=1/2$, it follows that each 
element of $b\in B_+$ has rank $\leq n$.
Thus,
$\Hi$ has dimension $k \leq n$,
and $B\cong M_k$. 
But the latter case was excluded in 
(iii,$n$).
Therefore $B$ is non-elementary 
and hence has the
\emph{global Glimm halving property}
in sense of 
\cite[def.~2.6]{BlanKir2}.
This is easy to see for
non-elementary simple $B$ (or  
use Glimm halving
\cite[lem.~6.7.1]{Ped.book}).
Since $B$ is simple property (iii, $n$) 
ensures that $B$ satisfies property 
(i) of \cite[def.~1.2]{BlanKir2} of pi($n$).
Therefore, \cite[prop.~4.14]{BlanKir2}
says that $B$ is pi($1$).  
Since $B$ is simple and $B\neq \mathbb{C}$, 
there are no non-zero characters on $B$.
In particular pi($1$) ensures that $B$ is 
purely infinite in the sense of 
\cite[def.~4.1]{KirRor1}. By 
\cite[thm.~4.1]{KirRor1} we obtain
property (ii).\qed

\smallskip

\textbf{Proof of Theorem 
\ref{thm:n.majorizing.n.covering}}:
By Lemma  
\ref{lem:deny.socle.for.major.cover}, 
$A$ is $G$-simple.
Thus the action $\sigma$ is 
automatically exact by 
Definition \ref{def:exact}.
Since $\sigma$ is  
element-wise properly outer by assumption, 
it is now
also residually properly outer, and
Theorem \ref{thm:A-filling} 
applies. It says that
$\cF:= A_+$ is a filling family in 
$A\rtimes_{\sigma,\lambda} G$.
In particular, $A$ 
separates the ideals of the 
reduced crossed product.
It shows that  
$B:=A\rtimes _{\sigma,\lambda} G$ is simple.

If $b\in B_+\,$ with $\| b\|=1$, 
then there exists non-zero
$z\in B$ such that $z^*z \leq b$ and 
$zz^*\in A$, 
because  $A_+$ is filling
for $B$:
One way to see this is to use 
property (i) of Definition
\ref{def:filling.mdiag.family}
on $a':=(3b-2)_+$, $b':=(3b-1)_+-a'$
and $c':=3b-a'-b'$ to get elements
$z_j, d\in B$. This imply that
$z_j^*z_j\leq e=\sum_i z_i^*z_i
=c'ec'\leq\|c'\|\|e\|c'
\leq 3\|c'\|\|e\|\,b$, where $e\neq 0$
because $d^*ed=a'\neq 0$.

Take  $\delta \in (0, \| z\|^2/2)$. 
Then
$0\not=z(z^*z-\delta)_+z^*=
\varphi(zz^*)\in A_+$
for some suitable 
$\varphi\in \Cf_c(0,\infty]_+$. 

We consider three cases: (i) the action is 
$n$-majorizing and $A$ is non-unital, (ii) 
the action is $n$-majorizing and $A$ is 
unital and (iii) the action is $n$-covering 
and $A$ is unital.

Since $G$ is discrete, 
$A$ is a non-degenerate 
\cst-subalgebra of $B$.
In particular, $A_+$
contains  an approximate unit $(e_\nu)$
of positive contractions in $A_+$ for $B$,
which we will use for case (i).
In case (ii)-(iii) let $e_\nu:=1$, where 
1 is the unit of 
$A$. Define $m:=n+1$ for case (ii) and 
$m:=n$ for case (i) and (iii).
By each of the 
Definitions 
\ref{def:n.majorizing.action} 
and 
\ref{def:n.covering.action}
(and using Lemma 
\ref{lem:deny.socle.for.major.cover}
to get property ($\beta$) in case (ii)
with $n$ replaced by $m$),
for each $\varepsilon >0$ and 
$e_\nu \in A_+$ 
there are
$d_1,\ldots,d_m\in A$ and 
$g_1,\ldots,g_m\in G$
such that 
$\| e_\nu - 
\sum_j f_j^*(z(z^*z-\delta)_+z^*)f_j \| 
< \varepsilon$
for $f_j:= U(g_j ^{-1})d_j$ in $B$.

By Remark \ref{rems:b-eps=d*ad}(ii), 
there exists a contraction
$d_0\in B$ with 
$d_0^*bd_0=(z^*z-\delta)_+$.
Then the elements 
$y_j:= d_0z^*f_j\in B$ satisfy
$\| e_\nu - \sum_{j=1}^{m}
y_j^* b y_j \| 
< \varepsilon$.
Since 
$a^{1/2} e_{\nu}a^{1/2} $ converges to 
$a\in B_+$,
we get  that $B$ has the  following property:

\emph{For any two non-zero elements 
$b,a \in B_+$
and $\varepsilon>0$ there exists 
$m$ elements
$d_1,\ldots, d_m\in B$ such that
$\|\, a\, - \,\sum _{j =1}^m  d_j^*bd_j\, \|\, 
<\,\varepsilon$.} A simple 
\cst-algebra $B$ with this property
is strongly purely infinite by 
Remark \ref{rem:simple.spi}
if $B$ is not isomorphic to $M_k$
for some $k\leq m$.
The case that $B$ is isomorphic to
a \cst-subalgebra of $M_m$ 
has been excluded by
the definitions: 
For case 
(i) we know $A$ is non-unital, hence $B$ 
is non-unital. For case (ii)-(iii) we know $A$ 
is not isomorphic to a \cst-subalgebra of $M_{m}$, 
hence $B$ is not 
isomorphic to a \cst-subalgebra of $M_{m}$.
\qed

\begin{lem}\label{lem:trivial.socle}
The following are equivalent for 
\cst-algebras $B$.
\begin{itemize}
\item[(i)] $B$
does not contain a non-zero projection
$p\in B_+$ with $pBp=\mathbb{C}\cdot p$.
\item[(ii)] For every 
non-zero hereditary \cst-subalgebra
$D$ of $B$, each
maximal commutative \cst-subalgebra $C$
of $D$ has perfect
Gelfand spectrum $\widehat{C}$.
\ie $\widehat{C}$ does not contain an isolated
point.
\item[(iii)]
For every $a_1,a_2\in B_+\backslash \{ 0\}$ and 
$c\in B$  there are 
$b_1,b_2\in B_+\backslash \{ 0\}$ with 
$b_1cb_2=0$,  and $b_j\leq a_j$  ($j=1,2$).
\end{itemize}
\end{lem}
\pf 
(iii)$\Rightarrow$(i):
Let $p^*=p=p^2\in B\backslash \{ 0\}$, put 
$a_k:=c:=p$ for $k=1,2$.  
Then there are non-zero
$b_1,b_2\in (pBp)_+$ with $b_1b_2=0$.
Thus $pBp \not= \mathbb{C}\cdot p$.

(i)$\Rightarrow$(ii):
Let $D\not=\{0\}$ a
hereditary \cst-subalgebra of $B$,
and $C$ a
maximal commutative
\cst-subalgebra of $D$. 
Suppose that $\widehat{C}$ is not perfect. 
Then
$\widehat{C}$
contains an isolated point $\chi$. 
The point $\chi$ corresponds
to a projection 
$0\not=p\in C_+\subseteq  D_+$ 
with $pBp=\mathbb{C}\cdot p$,
see \cite[lem.~7.14]{Ror.book}.

(ii)$\Rightarrow$(iii):
It is easy to see that a 
commutative \cst-algebra $C$
with a perfect  
spectrum $\widehat{C}$ 
contains non-zero
contractions $e_1,e_2\in C_+$ with $e_1e_2=0$,
because the locally compact Hausdorff space 
$\widehat{C}$ must in particular contain
two different points ($\not=\infty$).

Given $c\in B$ and non-zero $a_1,a_2\in B_+$, 
let
$d_j:= (a_j- \| a_j\|/2)_+$ and 
$x:=d_1^{1/2}cd_2^{1/2}$.
Notice that $0\not=d_j\in B_+$, 
and that $0\not= d_j^{1/2}y_jd_j^{1/2}\leq a_j$ 
for all
non-zero contractions
$0\leq y_j\in \overline{d_jBd_j}$.
If $x=0$ then take $b_j:=d_j$. 
If $x\not=0$, consider the hereditary
 \cst-subalgebra
 $D:=\overline{x^*Bx}=\overline{x^*xBx^*x}$
that is generated by $x^*x$, and is contained
in $\overline{d_2^*Bd_2}$.
Let  $C$ be a 
maximal commutative \cst-subalgebra
of $D$ with $x^*x\in C$.
Then $C$ contains non-zero contractions 
$e_1,e_2\in D_+$
with $e_1e_2=0=e_1(x^*x)^{1/2}e_2$. 

It is well-known (and easy to
see) that
the polar decomposition $x=v(x^*x)^{1/2}$
in $B^{**}$  has the property  that
$vDv^*=\overline{xBx^*}
\subseteq  \overline{d_1^*Bd_1}$. 
Thus
$f:=ve_1v^*\in \overline{xBx^*}$ 
and has the property 
$fxe_2=ve_1(x^*x)^{1/2}e_2=0$.
It follows that  $b_1:=d_1^{1/2}fd_1^{1/2}$
and 
$b_2:=d_2^{1/2}e_2d_2^{1/2}$
satisfy $b_1cb_2=d_1^{1/2}fxe_2d_2^{1/2}=0$
and $0\not=b_j\leq a_j$.\qed

\begin{prop}
\label{prop:1.majorizing.implies.G.separating}
Let $(A,G,\sigma)$ a \cst-dynamic system
with non-zero non-unital $A$. 
If the action
$\sigma$ of $G$ on $A$  is 
an $1$-majorizing  action 
in the sense of Definition 
\ref{def:n.majorizing.action},
then $A$ is $G$-simple and $\sigma$ is 
$G$-separating for $A$.
\end{prop}

\pf  The algebra $A$ is $G$-simple and 
$A$ does not contain
a projection $p\not=0$ with $pAp=\C\cdot p$
by Lemma \ref{lem:deny.socle.for.major.cover}.
To show 
that $\sigma$ is $G$-separating let
$a_1,a_2\in A_+\backslash  \{ 0\}$, $c\in A$ and 
$\varepsilon>0$. By Lemma \ref{lem:trivial.socle}, 
there exist $b_1,b_2\in A_+\backslash \{ 0\}$ with 
$b_1cb_2=0$ and $b_j\leq a_j$ ($j=1,2$). Using 
(twice) that the action is $1$-majorizing we find 
$e_j\in A$,  $h_j\in G$ for $j=1,2$ such that
$\| \, e_j^*\, \sigma_{h_j} \bigl(
b_ja_jb_j
\bigr) \,e_j\, - \, a_j\, \| \, 
<\, \varepsilon
\,$.
With $g_j:=h_j^{-1}$ and 
$d_j:=b_j\sigma_{g_j}(e_j)$ we get 
$\|d_j^*a_jd_j-\sigma_{g_j}(a_j)\|<\varepsilon$ 
and $d_1^*cd_2=0$, \ie  $\sigma$ 
is $G$-separating.\qed


\begin{rem}
\label{rem:unital.abel.2.cover=1.major=sba}
Suppose that $(A,G,\sigma)$ is a 
\cst-dynamical system, 
$A$ is unital and commutative, and 
$G$ is discrete.
Then the following 
properties (i)--(iv) of the action 
$\sigma$
are equivalent:
\begin{itemize}
\item[(i)] The action is 
$2$-covering in sense of 
Definition \ref{def:n.covering.action}.
\item[(ii)] The corresponding (adjoint)  action 
$\widehat{\sigma}$, on $\widehat{A}$ 
is a strong boundary action in the sense of 
Definition 
\ref{def:strong.boundary.action}.
\item[(iii)] The action is 
$1$-majorizing  in sense of 
Definition \ref{def:n.majorizing.action}.

\item[(iv)] The action is
$2$-filling in sense of 
Definition \ref{def:n.filling.action},
and $A$ is not isomorphic to a 
subalgebra of $M_{2}(\C)$.
\end{itemize}
We do not know if, also for 
non-commutative
and unital $A$, 
every $2$-covering action  is 
a  $1$-majorizing  action, 
or a $2$-filling action.%
\end{rem}
\pf We show more general implications, 
except for (i)$\Rightarrow$(ii).  In 
particular we show that an action 
$\sigma$ on a unital abelian 
\cst-algebra $A$ is $n$-filling if and 
only if it is $n$-covering 
provided that the (linear) dimension of 
$A$ is greater than $n$.

(i)$\Rightarrow$(iv) (for $A$ unital, 
commutative, any $n$):\\ 
Suppose that $A\cong \Cf (X)$, and take any 
$n\geq 2$. Let $\alpha$ denote the action of 
$G$ on $X$ inducing $\sigma$, i.e., 
$\sigma_g(f)=f\circ \alpha_g^{-1}$ for $f\in A$ 
and $g\in G$. Since the action 
$\alpha$ is minimal by Lemma 
\ref{lem:deny.socle.for.major.cover},
Remark \cite[rem.~0.4]{JolissaintRobertson}
shows it suffices to prove that
for each non-empty subset $U$ of $X$ 
there exist $g_1,\dots,
g_n\in G$, such that $\alpha_{g_1}(U)\cup 
\alpha_{g_2}(U) \cup \dots \cup 
\alpha_{g_n}(U)=X$.
Let such $U$ be given. Select non-zero
$a\in A_+$ with support contained in $U$.
By (i) there exist $g_1,\dots,g_n\in G$
and $d_1,\dots,d_n\in A$ such that 
$\sum_{j=1}^n d_j^*\sigma_{g_j}(a)d_j\geq 
\frac{1}{2}$. In particular for each
$x\in X$, $\sigma_{g_j}(a)(x)$ is non-zero
for some $j$, so $x\in \alpha_{g_i}(U)$.

(iv)$\Rightarrow$(i) (for $A$ unital, 
commutative/non-commutative, any $n$):\\ 
Suppose that $A$ is unital, and take any 
$n\geq 2$. Let $0\not=a\in A_+$.
Using (iv)
there are $g_1,\dots ,g_n\in G$
and $\delta > 0$ such that 
$\, D:= 
\sum_{j=1}^n\, \sigma_{g_j} (a)\, \ge\, \delta 1
\,,$
and $A$ is not isomorphic to a 
\cst-subalgebra of $M_n$. Thus,
 $D$ is invertible in $A$  and  
$\,
\sum_{j=1}^n \,  d_j^*\sigma_{g_j} (a) d_j\, = 1
\,$
for $d_j:=D^{-1/2}$ in $A$.

(iii)$\Rightarrow$(i) (for $A$ unital, 
commutative/non-commutative, any $n$):\\ 
Each $n$-majorizing actions 
on unital $A$ is an $(n+1)$-covering
action by 
Lemma
\ref{lem:deny.socle.for.major.cover}.

(i)$\Rightarrow$(ii) (for $A$ unital, 
commutative, one $n$):\\
Suppose that $A\cong \Cf (X)$ and let 
$\alpha$ denote the action of $G$ on $X$ 
inducing $\sigma$.
The equivalence of (i) and (iv),
shows that
for given $0\not=a\in A_+$, there exists
$g\in G$
and $\delta >0$ such that  
$a+\sigma_g (a)\ge \delta 1$. 

Let $U\subseteq X$ open and non-empty. 
There is 
$0\not= a\in \Cf (X)_+$ with
support $a^{-1}(0,\infty) \subseteq U$.
Choose $h\in G$ and $\delta>0$  with 
$a+\sigma_h (a)\ge \delta 1$.
It implies  that $\sigma_h (a)(x)>0$ 
for all $x\in X\backslash U$.
Thus, $\alpha_g(x)\in U$ for all 
$x\in X\backslash U$ and $g:=h^{-1}$,
\ie there exists $g\in G$ with 
$\alpha_g(X\backslash U)\subseteq U$.

Given non-empty
open  subsets $U$ and $V$ of $X$.
We let $W:= U\cap V$ if 
$U\cap V\not=\emptyset$.
Then $g\in G$ with 
$\alpha_g (X\backslash W)\subseteq W$
satisfies $\alpha_g (X\backslash U)\subseteq V$.
If $U\cap V=\emptyset$ then we find
$g,h \in G$ with 
$\alpha_g(X\backslash U)
\subseteq U\subseteq X\backslash V$
and $\alpha_h(X\backslash V)\subseteq V$.
Then 
$\alpha_{hg}(X\backslash U)\subseteq V$.

The space $X$ 
contains more than two points because 
$A$ is not isomorphic to a 
\cst-subalgebra of $M_2(\C)$.
Thus, $(X, G, \alpha)$ 
satisfies the conditions of 
Definition 
\ref{def:strong.boundary.action}
of a strong boundary action.

(ii)$\Rightarrow$(iii) (for $A$ 
unital/non-unital, commutative, one $n$):\\ 
We show (iii) using possibly less 
than (ii): Let
$X$ be a locally compact space
that is not necessarily compact 
and contains more than $2$ points.
Let $\alpha$ an action of $G$ on 
$X$ with the property that,
for every compact subset $K\subseteq X$
with $K\not=X$ and each 
non-empty open subset
$U\subseteq X$, there exists
$g\in G$ with $\alpha_g(K)\subseteq U$.
\emph{Then the adjoint action $\sigma$
of $\alpha$ 
on $A:=\Cf_0(X)$  is an $1$-majorizing
action of $G$  on $A$.}

Indeed: Let $0\not=a\in A_+$,  
$b\in A_+$ non-invertible, and 
$\varepsilon>0$.
Put $\delta:=\varepsilon/3$.
Then, considered as functions on $X$, 
they have the
property that $U:=a^{-1}( \| a\|/2, \infty)$ 
is non-empty and open and  
$K:= b^{-1}[\delta,\infty)$
is compact. 
Find $h\in G$ with 
$\alpha_{h}(K)\subseteq U$.
Then
$x\in K\,\,\Leftrightarrow\,\, b(x)\ge \delta$ 
implies that for $g:=h^{-1}$ we get
$\sigma_g (a)(x)=a(\alpha_{h}(x))>\| a\|/2\,$.
It follows 
$\| d^*\sigma_g(a) d -b \| <\varepsilon$  
with  $d\in A_+$
given by
$d(x):= 
\sigma_g(a)(x)^{-1/2} (b(x)-2\delta)_+^{1/2}$ 
for $x\in K$
and $d(x):=0$ for $x\in X\backslash K$.
\qed


\begin{rems}\label{rem:strong.boundary.action}	
(i) 
Let $\alpha$ be an action of a
discrete group $G$ on a locally compact Hausdorff
space $X$ with more than two points, and $\sigma$
the induced action on $A:=\Cf_0(X)$. By Remark
\ref{rem:unital.abel.2.cover=1.major=sba} the
following properties are equivalent for $X$ 
compact:  
\begin{itemize}
\item[(1)] The action $\alpha$ is a strong 
boundary action 
(Definition
\ref{def:strong.boundary.action})
in the sense of 
\cite{LacaSpi:purelyinf}: For each pair of 
non-empty open subsets $U$ and $V$ 
of $X$ there exists $g\in G$ with 
$U\cap \alpha_g(V)=X$.
\item[(2)] For any compact set $K\neq X$ and any 
open set $U\neq \emptyset$ there exist $t\in G$ 
such that $\alpha_t(K)\subseteq U$.
\item[(3)] For every non-zero $a\in A_+$, 
every non-invertible $b\in A_+\,$ and every 
$\varepsilon >0$, there exist
$d \in A$ and $g\in G$ 
such that
$\|d^*\sigma_{g} (a)d-b\|<\varepsilon$.
\end{itemize}
Clearly, this can not be the case if $X$ is 
locally compact but is not compact. In 
general 
(i.e., when $X$ is compact or non-compact)
we know 
(1)$\Rightarrow$(2)$\Rightarrow$(3). Properties
(2)-(3) are both candidates for a generalised 
notion of a strong boundary action, however 
only (3) applies when $A$ is non-commutative.

\noindent
(ii)
The notion of a strong 
boundary action 
(Definition \ref{def:strong.boundary.action}) is
defined on \emph{compact} Hausdorff spaces 
with more than two points. In view of Remark
\ref{rem:strong.boundary.action}(i) and Remark
\ref{rem:unital.abel.2.cover=1.major=sba}, 
we call the $1$-majorizing actions on 
not-necessarily unital or commutative 
\cst-algebras also strong boundary 
actions.

\noindent
(iii) 
Suppose that a discrete 
group $G$ acts by a topologically free 
action $\alpha$ on a \emph{compact} Hausdorff 
space $X$, and that $X$ contains more than two 
points. It was shown in 
\cite[thm.~5]{LacaSpi:purelyinf} that 
the crossed product 
$\Cf (X)\rtimes_{\sigma,\lambda} G$ 
is purely infinite 
provided that the action -- in addition -- is a 
strong boundary action. Since topological 
freeness implies $\sigma$ is element-wise properly 
outer (by \cite[prop.~1]{ArchSpiel}) we conclude 
that, with the terminology of Remark 
\ref{rem:strong.boundary.action}(ii), 
\cite[thm.~5]{LacaSpi:purelyinf} is a special case 
of Theorem \ref{thm:n.majorizing.n.covering}.
\qed

\noindent
(iv) 
Let $\alpha$ be an action on 
a \emph{non-compact} locally compact Hausdorff space $X$ 
with more than two points and $\sigma$ the induced 
action.
It was shown in Proposition 
\ref{prop:1.majorizing.implies.G.separating} that 
$\sigma$ is $G$-separating if $\sigma$ is a strong 
boundary (\ie $1$-majorizing) action. A simpler 
argument applies if we assume that for any compact 
set $K\neq X$ and any non-empty open set $U\subseteq X$ 
there exist $g\in G$ such that 
$\alpha_g(K)\subseteq U$:\\
\pf
Since any finite subset $M$  of 
$X$ is compact, it can be
moved by suitable 
$\alpha_g$ into any non-empty open subset 
$U$ of $X$.  
In particular $X$ is perfect and each non-empty
open set
$V\subseteq X$ contains at least two
non-empty open disjoint subsets 
$V_1$ and $V_2$.
Let $K_1\subseteq U_1$ and $K_2\subseteq U_2$  
with $K_j$ compact (hence $K_j\neq X$)
and $U_j$ open. 
If $U_1$ and $U_2$ are disjoint, then
we can take $g_1=g_2=e$ in  Lemma
\ref{lem:G-separating}(ii).
If $V:= U_1\cap U_2\not=\emptyset$, then 
consider the above disjoint 
non-empty open subsets 
$V_j\subseteq V$.
By assumption, there exist $g_1,g_2\in G$ with
$\alpha_{g_j}(K_j)\subseteq V_j \subseteq U_j$.
Thus, the adjoint action $\sigma$
of $\alpha$ is $G$-separating.
\qed

\noindent
(v) 
Suppose that a discrete 
group $G$ acts by a topologically free 
action $\alpha$ on a \emph{non-compact} 
locally compact Hausdorff space $X$, 
and that $X$ contains more 
than two points. Then the crossed product 
$\Cf_0(X)\rtimes_{\sigma,\lambda}  G$ is purely 
infinite provided that the following 
property holds: for any compact set $K\neq X$ and 
any non-empty open set $U\subseteq X$ there exist 
$t\in G$ such that $\alpha_t(K)\subseteq U$. 
This follows as a corollary of Theorem 
\ref{thm:main}, also of Theorem 
\ref{thm:n.majorizing.n.covering} or of 
Corollary \ref{cor:Abelian}.\\
\pf
We must verify the following properties according 
to each of the listed results:
\begin{itemize}
\item[({\ref{thm:main}})] The action $\sigma$ is 
exact, residually properly outer and 
$G$-separating.
\item[({\ref{thm:n.majorizing.n.covering}})] The 
action $\sigma$ is $1$-majorizing, and 
element-wise properly outer.
\item[(\ref{cor:Abelian})] The action $\sigma$ 
is exact, $G$-separating and fulfills (*): For 
every closed $G$-invariant subset $Y$ of $X$ and 
every $g\neq e$
the set $\{y \in Y\colon 
\alpha_g(y)=y\}$ has empty interior.
\end{itemize}
By Remark \ref{rem:strong.boundary.action}(i) we 
know the action $\sigma$ is a strong boundary 
(\ie $1$-majorizing) action. Hence $A$ is 
$G$-simple, \cf{}Lemma 
\ref{lem:deny.socle.for.major.cover}. In 
particular it follows that  the action $\alpha$ 
on $X$ is minimal. This reduces  property 
(*) to the definition of topological freeness, 
\cf{}\cite[p.120]{ArchSpiel}. The minimality of
the action $\alpha$ implies that the corresponding
\emph{adjoint action} 
$\sigma\colon G\to \mathrm{Aut}(\Cf_0(X))$ 
is exact,
and that it becomes  residually properly outer 
if it is element-wise properly outer. But $\sigma$ 
is element-wise properly outer if and only if 
$\alpha$ is a topological free action 
(see \cite[p.120]{ArchSpiel} or 
\cite[cor~2.22]{Siera2010}). It remains to show 
$\sigma$ is $G$-separating, but this is already 
contained in Remark 
\ref{rem:strong.boundary.action}(iv).
\qed

\noindent
(vii)
It is an important point that a strong 
boundary action is often $G$-separating 
and (in fact always) minimal, but the 
notion of a $G$-separating action is not typically 
related to minimality. Consequently, working with 
$G$-separating 
actions
allows us to consider ideal-related classification 
of non-simple strongly purely infinite 
crossed products.
\end{rems}

\begin{rem}\label{rem:case.RR=0}
If $A$ has real rank zero, then one can restrict 
the conditions in Definitions
\ref{def:G.separating},  
\ref{def:n.majorizing.action} and 
\ref{def:n.covering.action}
to projections $p,q\in A$ in place of the elements  
$a,b\in A_+$.
\end{rem}

\pf  Case of Definition \ref{def:G.separating}:\,
Let $a_1, a_2\in A_+$,  $c\in A$,  
$\varepsilon>0$ 
and define
$$\delta:=\varepsilon/(1+\| a_1\| + \| a_2\|)\,.$$
By \cite{BrownPedersen:RR=0}, 
$D_j := \overline{\, a_jAa_j\, }$ 
contains 
an approximate unit consisting of non-zero 
projections. Thus, there are
projections $p_j\in D_j$ such that
$\| a_j-\, a_j^{1/2}\, p_j\, a_j^{1/2} \|<\delta$.
Use \cite[prop.~2.7(i)]{KirRor1} and the comment
following \cite[prop.~2.6]{KirRor1} to
select $z_j\in D_j$
satisfying $z_j^*a_j\,z_j\,=p_j $. 
Let $ c':=z_1^*cz_2 $. 

Suppose that
there exists $e_j\in A$, $g_j\in G$ such that
$$\| e_j^*p_j e_j  - \,\sigma_{g_j} (p_j)\| < 
\delta \quad \text{and}
\quad \| e_1^*c'\, e_2\| <\delta \,.$$
Define $v_j := \sigma_{g_j}(a_j^{1/2})$ and
$d_j := z_j e_j \sigma_{g_j} (a_j^{1/2})\, $.  
They satisfy
$d_j^*a_jd_j= v_j^* e_j^*p_je_j v_j$,  
$v_j^*\sigma_{g_j}(p_j) v_j = 
\sigma_{g_j}(a_j^{1/2} p_j a_j^{1/2})$.
Thus, 
$
\| d_j^*a_j d_j - \sigma_{g_j}(a_j) \| 
< (1+ \| a_j\|) \delta \leq  \varepsilon 
\,$.
Since 
$d_1^*cd_2 =    v_1^* e_1^* c' e_2 v_2$,  
we get
$\|  d_1^* c d_2  \| 
< \delta (\| a_1\| \cdot  \| a_2\|)^{1/2}\leq 
\delta (\| a_1\| +  \| a_2\|)\leq
\varepsilon\,$.

Case of Definitions \ref{def:n.majorizing.action}  
and \ref{def:n.covering.action}:\,\, 
Let $a_1, a_2 \in A_+$ and 
$\varepsilon >0$, with $a_1\not=0$ and
$a_2$ not invertible in 
$A$ (respectively $a_2=1$).
We can assume $\varepsilon \leq 1$.
Define
$\delta := \varepsilon /(1+ \| a_2\|)\,$.
Choose 
$p_j,z_j \in D_j:= \overline{ a_jAa_j}$ 
as above with
$\| a_j^{1/2} p_j a_j^{1/2}- a_j\| <\delta$ and 
$z_j^*a_jz_j=p_j$.
Then 
$p_1\not=0$ and $p_2$ is not invertible in 
$A$ (\ie $p_2\not=1$)
if $a_2$ is not invertible, otherwise $p_2=1$:
If $p_2$ is invertible then
$1\in a_2Aa_2$, so $a_2$ is invertible. Conversely, 
if $a_2$ is invertible then 
$\|p_2-1\|<\varepsilon/2$, so $p_2$ is invertible.

If there are 
$e_1,\ldots, e_n \in A$ and $g_1,\ldots, g_n\in G$ 
with
$\| \, p_2\, -\, 
\sum_j \, e_j^* \sigma_{g_j}(p_1)e_j \, \|
 < \delta\,,
$
then 
$d_j :=\sigma_{g_j}(z_1) e_ja_2^{1/2}$ 
satisfies 
$
\| \, a_2\, -\, \sum_j \, d_j^* 
\sigma_{g_j}(a_1)d_j \, \|  
< 
(1+ \|a_2\|)\delta \leq \varepsilon \,$.
\qed

\section*{Acknowledgments}
Parts of this work were conducted while the second 
named author was at the Fields Institute from 2009
to 2012. It is with great pleasure we forward our 
thanks to the Fields Institute and in particular 
Professor George Elliott for all the support.
This research was also supported by the 
Australian Research Council.

\appendix
\section{}
\label{appendix}
In this appendix we have included a few recent 
definitions and results that are frequently cited 
throughout this paper. The results quoted from 
{\cite{KirSie2}} are available as preprint.  

\begin{definition}[{\cite{LacaSpi:purelyinf}}]
\label{def:strong.boundary.action}
Let $\alpha$ be an action of a discrete group
$G$ on a compact spaces $X$ with at least three 
points. The action $\alpha$ is as
\emph{strong boundary action} if for every pair
$U$ and $V$ of non-empty open subsets of $X$ there
exists $t\in G$ such that 
$\alpha_t(X\setminus U)\subseteq V$.
\end{definition}

\begin{definition}[{\cite{JolissaintRobertson}}]
\label{def:n.filling.action}
An action $\sigma$  of a discrete group $G$ on 
a unital \cst-algebra $A$ is \emph{$n$-filling} 
($n\geq 2$) if, for all $b_1,\dots, b_n\in A_+$,
with $\|b_j\|=1$ for each $j$, and all $\varepsilon>0$,
there exist $g_1,\dots,g_n\in G$ such that
$\sum_{j=1}^n \sigma_{g_j}(b_j)\geq 1-\varepsilon$.
\end{definition}	

\begin{definition}[{\cite{KirSie2}}]
\label{def:filling.mdiag.family}
Let $\cF $  be a subset of  $A_+$.
The set $\cF $ is a 
\emph{filling family} for $A$, if
$\cF $ satisfies the following equivalent 
conditions (i) 
and (ii).
\begin{itemize}
\item[(i)]
For every $a,b,c\in A$ with 
$0\leq a\leq b\leq c\leq 1$,
with 
$ab=a\not=0$ and  $bc=b$, there exists
$z_1,z_2,\ldots, z_n\in A$ and $d\in A$ with
$z_j(z_j)^*\in \cF $, such that
$ec=e$
and $d^*ed=a$ for 
$e:=z_1^*z_1+\ldots+z_n^*z_n$.
\item[(ii)] 
For every hereditary \cst-subalgebra
$D$ of $A$ and every primitive ideal
$I$ of $A$ with $D\not\subseteq I$ there
exist $f\in \cF $ and $z\in A$
with $z^*z\in D$ and  $zz^*=f  \not\in I$.
\end{itemize}
\end{definition}

\begin{lem}[{\cite{KirSie2}}]
\label{lem:F.fill.A.A+.fill.B.Then.F.fill.B}
Suppose that $A\subseteq B$ is a 
\cst-subalgebra of $B$
and
$\cF\subseteq  A_+$ is a subset of $A_+$.
If $\cF$ is filling for 
$A$, and $A_+$ is filling for $B$, then $\cF$
is a filling family for $B$.
\end{lem}

\begin{rem}[{\cite{KirSie2}}]
\label{rems:filing.familiy}
Let $A\subseteq B$ be \cst-algebras
and $\cF \subseteq  A_+$.
If
$\cF:=A_+\subseteq B$ is filling for $B$, then
the map 
$I\in \mathcal{I}(B)\mapsto
 I\cap A \in  \mathcal{I}(A)$ 
is injective, \ie $A$ separates the 
closed ideals
of $B$.
\end{rem}

\begin{definition}[{\cite{KirSie2}}]
\label{def:spi}
A \cst-algebra $A$ is
\emph{strongly purely infinite}
(for short: \textit{s.p.i.}\,)  if,
for every $a, b\in A_+\,$ and  $\varepsilon>0$, 
there exist elements $s,t\in A$ such that
\begin{equation}
\label{InEq.spi}
\| \,  s^*a^2s - a^2 \,\|  \, <\varepsilon\,,\;
\| \,  t^*b^2t - b^2\,\|  \, <\varepsilon
\;\;\textrm{and}\;\; \| \, s^*abt\,\|  \, 
<\varepsilon\, .
\end{equation}
\end{definition}

\begin{rem}[{\cite{KirSie2}}]
\label{rem:spi}
A \cst-algebra $A$ is
strongly purely infinite if and only if
for every $a, b\in A_+$, $c\in A$ and  $\varepsilon>0$, 
there exist contractions $s,t\in A$ such that
\begin{equation}
\label{InEq.spi2}
\| \,  s^*as - a \,\|  \, <\varepsilon\,,\;
\| \,  t^*bt - b\,\|  \, <\varepsilon
\;\;\textrm{and}\;\; \| \, s^*ct\,\|  \, 
<\varepsilon\, .
\end{equation}
\end{rem}

\begin{definition}[{\cite{KirSie2}}]
\label{def:MatrixDiag}
Let $\cS\subseteq A$ 
be a multiplicative sub-semigroup 
of a \cst-algebra $A$ and $\Co \subseteq A$
a subset of $A$.
An $n$-tuple $(a_1,\ldots,a_n)$ of
positive elements in $A$ has the 
\emph{matrix diagonalization property 
with respect to $\cS$ and $\Co $},
if for every 
$[a_{ij}]\in M_n(A)_+$ with $a_{jj}=a_j$ and 
$a_{ij}\in \Co $
(for $i\not=j$)
and $\varepsilon_j >0, \tau>0$ 
there are elements 
$s_1\,,\ldots , s_n\, \in \cS$ 
with
\begin{equation}
\label{InEq.DiagProp}
\| s_j ^* a_{jj}\, s_j - a_{jj} \| 
< \varepsilon_j\,, \quad
\text{and} \quad 
\| s_i ^*a_{i j}s_j \| < \tau \,\,\,  
\text{for} \,\, i\not= j\,.
\end{equation}
If $\cS=\Co =A$ then this is
the \emph{matrix diagonalization} property 
of $(a_1,\ldots,a_n)$ as defined in 
\cite[def.~5.5]{KirRorOinf},
and we say that $(a_1,\ldots,a_n)$ 
has matrix diagonalization (in $A$).
\end{definition}

\begin{definition}[{\cite{KirSie2}}]
\label{def:mdiag.family}
Let $\cF $  be a subset of  $A_+$.
The family $\cF $ has the 
\emph{(matrix) diagonalization property}
(in $A$)  if each finite sequence 
$a_1,\ldots, a_n\in \cF $
has the matrix diagonalization property 
(in $A$) of 
Definition \ref{def:MatrixDiag}.
\end{definition}

\begin{lem}[{\cite{KirSie2}}]
\label{lem:flexible.diag.imply.general}
Suppose that $\cF\subseteq  A_+$ is invariant
under $\varepsilon$-cut-downs, \ie that
for each $a\in \cF$ and 
$\varepsilon \in (0,\| a\|)$
we have
$(a-\varepsilon)_+\in \cF$.

Then the family  
$\cF$ has the matrix diagonalization property, 
if and only if, each pair 
of elements in $\cF$  
has the matrix diagonalization property of 
Definition \ref{def:MatrixDiag}.
\end{lem}

\begin{lem}[{\cite{KirSie2}}]
\label{lem:cS.control}
Let $\varepsilon_0>0$
and non-empty subsets $\cF\subseteq A_+$,
$\Co\subseteq A$
be given, and let 
$\cS \subseteq A$ be a 
(multiplicative) sub-semigroup of $A$ 
that satisfies $s_2^*\Co s_1\subseteq  
\Co$ for all $s_1,s_2\in \cS$.
Suppose that
the following
properties hold:
\begin{itemize}
\item[(i)] For every 
$\varepsilon_0>\delta>0$, the pair 
$((a_1-\delta)_+, (a_2-\delta)_+)$  
 the matrix 
diagonalization
property with respect to $\cS$ and $\Co $
 of Definition 
\ref{def:MatrixDiag}.
\item[(ii)]
$\varphi(a_1)c\varphi (a_2)\in \Co$ 
for each $c\in \Co$
and $\varphi \in \Cf_c (0,\infty]_+$.
\item[(iii)]
$\varphi(a_1)s, \varphi(a_2)s\in \cS$ 
for each $s\in \cS$
and $\varphi \in \Cf_c (0,\infty]_+$.
\end{itemize}
Then, for every 
$c\in \overline{\mathrm{span}(\Co )}$,
$a_1,a_2\in \cF$,
$\varepsilon_0/2\geq\varepsilon >0$, and
$\tau>0$, 
there exists 
$s_1,s_2\in \cS$
that fulfil
$\|s_j\|^2 \leq 2\| a_j \|/\varepsilon$
and
\begin{equation}\label{InEq.general}
\| s_1^* a_1 s_1- a_1 \|<\varepsilon\,, \quad
\| s_2^* a_2 s_2- a_2 \|<\varepsilon\, \quad
\text{and} \quad \| s_1^*cs_2 \| <\tau\,.
\end{equation}
\end{lem}

\begin{thm}[{\cite{KirSie2}}]
\label{thm:A.ot.B.spi.if.A.spi.B.exact}
The minimal tenor product of a strongly purely 
infinite and an exact \cst-algebra is strongly 
purely infinite.
\end{thm}

\begin{thm}[{\cite{KirSie2}}]
\label{prop:local-spi}
Suppose that $A_+$ contains a filling family
$\cF$ that has the diagonalization property 
(in $A$).
Then $A$ is strongly purely infinite.
\end{thm}

\providecommand{\bysame}{\leavevmode\hbox
to3em{\hrulefill}\thinspace}

\bigskip

\address{Institut f{\"u}r Mathematik,
Humboldt Universit{\"a}t zu Berlin, 
Unter den Linden 6,
D--10099 Berlin, Germany}\\ 
\email{kirchbrg@mathematik.hu-berlin.de}

\address{School of Mathematics \& 
Applied Statistics,
University of Wollongong, 
Faculty of Engineering \& Information Sciences,
2522  Wollongong, Australia}\\
\email{asierako@uow.edu.au}

\begin{thebibliography}{10}
%
\bibitem{ArchSpiel}
R.J.~Archbold and J.S.~Spielberg, 
{\em Topologically free actions 
and ideals in discrete C*-dynamical systems}, 
Proc.~Edinburgh Math.~Soc.~(2) 
\textbf{37} (1994), 119--124.

\bibitem{Arzhantseva.Guba.Sapir}
G.~Arzhantseva, V.~Guba, and M.~Sapir,
{\em Metrics on diagram groups and uniform
embeddings in Hilbert space},
Comment.\ Math.~Helv. \textbf{81} (2006), 
911--929.

\bibitem{BlanKir2}
E.~Blanchard and E.~Kirchberg, 
{\em Non-simple purely infinite
C*-algebras: the Hausdorff case}, 
J.\  Funct.~Anal.\  \textbf{207} (2004), 
461--513. 

\bibitem{BrownPedersen:RR=0}
L.G.~Brown and G.K.~Pedersen,
\emph{C*-algebras of real rank zero},
J.\  Funct.~Anal.\  \textbf{99} (1991), 
131--149.

\bibitem{Cuntz:K-Th.in.Annals}
J.~Cuntz,
\emph{K-theory for certain C*-algebras},
Ann.~of Math. \textbf{113} (1981), 181--197.

\bibitem{Elliott.prop.outer}
G.A.~Elliott, 
\emph{Some simple 
C*-algebras constructed 
as crossed products with discrete outer
automorphism groups}, 
Publ. Res.~Inst.~Math.~Sci. \textbf{16} (1980), 
299--311.

\bibitem{Exel}
R.~Exel, M.~Laca, and J.~Quigg, 
\emph{Partial dynamical systems and
C*-algebras generated by partial isometries}, 
J.\ Operator Theory	 \textbf{47} (2002), 
169--186. 

\bibitem{FarleyDS.Thompson}
D.S.~Farley, 
\emph{Proper isometric actions of 
Thompson's groups on Hilbert space},
Int. Math. Res. Notes \textbf{45} (2003), 
2409--2414.

\bibitem{Guentner.Kaminker}
E.~Guentner and J.~Kaminker, 
\emph{Exactness and uniform 
embeddability of discrete Groups},
J.~London Math.~Soc.(2) 
\textbf{70} (2004), 
703--718.

\bibitem{HaagerupPicioroaga}
U.~Haagerup and G.~Picioroaga,
\emph{New presentations of 
Thompson's groups and applications},
J. Operator Theory \textbf{66} (2011), 
217--232.

\bibitem{JolissaintRobertson}
P.~Jolissaint and G.~Robertson,
\emph{Simple purely infinite 
C*-algebras and $n$-filling actions},
J. Funct. Anal. \textbf{175} (2000), 
197--213.

\bibitem{KawamuraTomiyama} 
S.~Kawamura and J.~Tomiyama, 
\emph{Properties of 
topological dynamical systems and 
corresponding C*-algebras}, 
Tokyo J.~Math. \textbf{13} (1990), 251--257.

\bibitem{KirRor1}
E.~Kirchberg and M.~R{\o}rdam, 
\emph{Non-simple purely infinite C*-algebras},
Amer.\  J.~Math. \textbf{122} (2000), 637--666.

\bibitem{KirRorOinf}
\bysame,
\emph{Infinite non-simple C*-algebras: absorbing 
the Cuntz algebras $\mathcal{O}_{\infty}$}, 
Adv.~Math. \textbf{167} (2002), no.~2, 195--264. 

\bibitem{KirSie2} 
E.~Kirchberg and A.~Sierakowski, 
\emph{Filling families and strong 
pure infiniteness}, preprint 2014.

\bibitem{KirWa.Exact.Groups.Cont.Bund}	  
E.~Kirchberg and S.~Wassermann, 
\emph{Exact groups and 
continuous bundles of C*-algebras}, 
Math.~Ann. {\bf 315} (1999), 169--203.

\bibitem{Kishimoto.outer.auto} 
A.~Kishimoto, 
\emph{Outer automorphisms 
and reduced crossed products of simple 
C*-algebras}, 
Comm.~Math.~Phys. 
\textbf{81} (1981), 429--435.

\bibitem{LacaSpi:purelyinf}
M.~Laca and J.~Spielberg, 
\emph{Purely infinite C*-algebras from
	  boundary actions of discrete groups}, 
J.\ Reine Angew.~Math. \textbf{480}  (1996), 
125--139. 

\bibitem{OlePed3}
D.~Olesen and G.K.~Pedersen, 
{\em Applications of the Connes spectrum
to C*-dynamical systems. \textrm{III}}, 
J. Funct. Anal. \textbf{45} (1982), 357--390. 

\bibitem{Ped.book}
G.K.~Pedersen, 
\emph{C*-algebras and their automorphism
groups}, 
LMS Monographs, vol.~14, 
Academic Press Inc., London, 1979. 

\bibitem{Ror.book}
M.~R{\o}rdam, F.~Larsen, and N.~Laustsen,
\emph{An introduction to K-theory 
for C*-algebras},
London Mathematical Society Student Texts
\textbf{49},
Cambridge University Press, Cambridge, 2000.

\bibitem{Ren:fixed}
J.~Renault, 
\emph{The ideal structure of 
groupoid crossed product 
C*-algebras},
\emph{With an appendix by 
Georges Skandalis.}
J. Operator Theory \textbf{25} (1991), 3--36. 

\bibitem{Siera2010}
A.~Sierakowski, 
\emph{The ideal structure of 
reduced crossed products}, 
M\"unster J.~Math. 
\textbf{3} (2010), 223--248.

%
\end{thebibliography}
\end{document}